\documentclass[11pt,letterpaper]{article}
\usepackage{fullpage}
\usepackage[margin=1.038in]{geometry}
\usepackage{multirow}
\usepackage{url}
\usepackage[utf8]{inputenc}
\usepackage{amsmath,amssymb,amsfonts}
\usepackage{graphicx}
\usepackage{color}
\usepackage{fancyhdr}
\usepackage{tikz}
\usepackage{soul}
\usepackage{comment}
\usepackage{bm}
\usepackage{enumitem,kantlipsum}
\usepackage{wrapfig}
\usepackage[compact]{titlesec}
\titlespacing*{\section} {0pt}{1ex}{1ex}
\titlespacing*{\subsection} {0pt}{1ex}{1ex}
\titlespacing*{\subsubsection}{0pt}{1ex}{1ex}

\usepackage{amsopn}

\usepackage{amsthm}
\usepackage[square,sort,comma,numbers]{natbib}
\DeclareMathOperator*{\argmin}{\arg\!\min}
\DeclareMathOperator*{\argmax}{\arg\!\max}

\renewcommand*{~}{\relax\ifmmode\sim\else\nobreakspace{}\fi}

\newcommand{\Fbar}{\bar{F}}

\newcommand{\Nhat}{\hat{N}}

\newcommand{\Rhat}{\hat{R}}

\newcommand{\sigmahat}{\hat{\sigma}}

\newcommand{\Ntilde}{\widetilde{N}}

\newcommand{\Rtilde}{\widetilde{R}}

\newcommand{\BFe}{\bm{e}}

\newcommand{\BFs}{\bm{s}}

\newcommand{\BFx}{\bm{x}}

\newcommand{\BFG}{\bm{G}}

\newcommand{\BFU}{\bm{U}}
\newcommand{\BFX}{\bm{X}}

\newcommand{\BFY}{\bm{Y}}

\newcommand{\BFbeta}{\bm{\beta}}

\newcommand{\BFXhat}{\hat{\BFX}}

\newcommand{\BFXtilde}{\widetilde{\BFX}}

\newcommand{\sfH}{\mathsf{H}}
\newcommand{\sfM}{\mathsf{M}}

\newcommand{\mcB}{\mathcal{B}}

\newcommand{\mcD}{\mathcal{D}}

\newcommand{\mcF}{\mathcal{F}}
\newcommand{\mcG}{\mathcal{G}}

\newcommand{\mcN}{\mathcal{N}}

\newcommand{\mcR}{\mathcal{R}}

\newcommand{\mcV}{\mathcal{V}}

\newcommand{\mcX}{\mathcal{X}}

\newcommand{\mbE}{\mathbb{E}}

\newcommand{\mbP}{\mathbb{P}}

\newcommand{\mbR}{\mathbb{R}}

\usepackage{array, caption, floatrow, tabularx, makecell, booktabs}%
\setcellgapes{3pt}
\usepackage{algorithmic}
\usepackage{algorithm}

\newtheorem{theorem}{Theorem}

\newtheorem{lemma}{Lemma}

\newtheorem{proposition}{Proposition}

\newtheorem{corollary}[theorem]{Corollary}

\newtheorem{definition}{Definition}

\newtheorem{assumption}{Assumption}

\newtheorem{remark}{Remark}

\newcommand{\TD}{\nabla}

\usepackage{bm}
\usepackage{multicol}
\usepackage{lipsum}
\usepackage{mwe}
\usepackage{graphicx}
\usepackage{subfig}
\usepackage{amssymb}
\usepackage{booktabs}
\usepackage{enumitem}
\usepackage{cleveref}

\usepackage{multirow}
\usepackage{amsfonts}
\usepackage{graphicx}
\usepackage{epstopdf}
\usepackage{algorithmic}
\usepackage{algorithm}

\usepackage{qcircuit}
\usepackage{physics}

\usepackage{authblk}

\ifpdf
  \DeclareGraphicsExtensions{.eps,.pdf,.png,.jpg}
\else
  \DeclareGraphicsExtensions{.eps}
\fi

\title{Two-Stage Estimation and Variance Modeling for Latency-Constrained Variational Quantum Algorithms}
\author[1]{Yunsoo Ha, Sara Shashaani}
\author[2]{Matt Menickelly}
\affil[1]{Department of Industrial and Systems Engineering, North Carolina State University}
\affil[2]{Mathematics and Computer Science Division, Argonne National Laboratory}
\date{}

\begin{document}
\vspace{-3 mm}

\maketitle

\begin{abstract}
The Quantum Approximate Optimization Algorithm (QAOA) has enjoyed increasing attention in noisy intermediate-scale quantum computing due to its application to combinatorial optimization problems. Because combinatorial optimization problems are NP-hard, QAOA could serve as a potential demonstration of quantum advantage in the future. 
As a hybrid quantum-classical algorithm, the classical component of QAOA resembles a simulation optimization problem, in which the simulation outcomes are attainable only through the quantum computer. 
The simulation that derives from QAOA exhibits two unique features that can have a substantial impact on the optimization process: (i) the variance of the stochastic objective values typically decreases in proportion to the optimality gap, and (ii) querying samples from a quantum computer introduces an additional latency overhead.
In this paper, we introduce a novel stochastic trust-region method, derived from a derivative-free adaptive sampling trust-region optimization (ASTRO-DF) method, intended to efficiently solve the classical optimization problem in QAOA, by explicitly taking into account the two mentioned characteristics. 
The key idea behind the proposed algorithm involves constructing two separate local models in each iteration: a model of the objective function, and a model of the variance of the objective function. 
Exploiting the variance model allows us to both restrict the number of communications with the quantum computer, and also helps navigate the nonconvex objective landscapes typical in the QAOA optimization problems. 
We numerically demonstrate the superiority of our proposed algorithm using the SimOpt library and Qiskit, when we consider a metric of computational burden that explicitly accounts for communication costs.
\end{abstract}

\section{Introduction}
Quantum computers have the potential to outperform their classical counterparts on numerous critical calculations.
Diverse fields including data science~\citep{biamonte2017quantum}, quantum chemistry~\citep{lanyon2010towards}, condensed matter~\citep{smith2019simulating}, nuclear physics~\citep{cloet2019opportunities}, and even finance~\citep{orus2019quantum} stand to benefit from quantum algorithms in various ways in the future. 
However, in the near-term, during the noisy intermediate-scale quantum (NISQ) era~\citep{preskill2018quantum}, realizing these theoretical advantages is challenging.
This is because canonical quantum algorithms used in many of these fields
necessitate gate depths that are only expected to be achievable with fault-tolerant, error-corrected quantum computers~\citep{preskill1998fault}.

Variational quantum algorithms (VQAs) aim to reduce gate depth requirements by exploiting classical computer-based optimization processes~\citep{cerezo2021variational}. 
These algorithms have demonstrated their effectiveness on NISQ hardware in tasks such as dynamical evolution~\citep{yuan2019theory,otten2019noise,kandala2017hardware}, eigenvalue estimation~\citep{o2016scalable}, machine learning~\citep{mitarai2018quantum,otten2020quantum}, and various other problem domains~\citep{cerezo2021variational}.
One of the primary challenges in VQAs lies in the optimization step, which is performed on classical computers. 
The optimization step involves estimating an expectation cost function  (and potentially, its derivative information) derived from the problem being solved, using a limited number of samples. 
These samples are often referred to as \emph{shots} in this context. 
The estimation of an expectation objective function ideally necessitates the employment of stochastic optimization algorithms.
A straightforward way to quantify the overall cost of a classical optimizer is by counting the total number of shots executed on the quantum device to 
achieve an $\epsilon$ zeroth-order or first-order optimality gap. 
To estimate the cost function with a given a set of parameters (decision variables), multiple shots must be executed on a correspondingly parameterized quantum circuit. The estimation error is quantifiable analogous to Monte Carlo estimators. 
With this perspective, the optimization performed on the classical computer can be seen as a form of simulation optimization (SO). 
For flexibility, we make no assumptions about the accessibility of (directional) derivatives in the VQA context; such a setting necessitates \emph{derivative-free} SO solvers. Derivative-free SO solvers generate solution paths for simulations (stochastic oracles) that do not provide direct derivative observations, also known as zeroth-order oracles. 
Before discussing SO in further detail, we begin by introducing a specific example of VQA especially relevant to operations research. 
\subsection{Quantum Approximate Optimization Algorithm}\label{sec:qaoa}
The quantum approximate optimization algorithm (QAOA) is a particular and well-studied instance of VQA, designed for the solution of a class of combinatorial optimization algorithms.
In QAOA, once a combinatorial optimization problem is fixed, a matrix $H_C$, called the cost Hamiltonian, is specifically (and implicitly) constructed in such a way to ensure that its ground state (lowest eigenvalue)  corresponds to the optimal solution to the original combinatorial optimization problem. 
QAOA relies on what is known as the ``variational principle," which states $\bra{\psi(\BFx)}H_C\ket{\psi(\BFx)} \ge E_0,$
where $E_0$ is the ground state energy and $\ket{\psi(\BFx)}$ is a quantum state vector parameterized by $\BFx$. We thus aim to solve the problem of the form
\begin{equation} \label{eq:qaoa-problem}
    \min_{\BFx\in\mbR^d}\bra{\psi(\BFx)}H_C\ket{\psi(\BFx)}.
\end{equation}
Given that the quantum state vector collapses to a single state upon measurement, we must estimate $\bra{\psi(\BFx)} H_C \ket{\psi(\BFx)}$, which represents an expectation of a physical quantity, by repeatedly measuring the quantum state and employing Monte Carlo sampling. 
For convenience of notation in discussing an optimization algorithm, we let $F(\BFx,\xi)$ denote a stochastic function value (all sources of stochasticity are encoded in the random variable $\xi$), and we can then rewrite \eqref{eq:qaoa-problem} as 

$$\displaystyle\min_{\BFx\in\mbR^d} f(\BFx) := \mbE_{\xi}[F(\BFx,\xi)].$$

\begin{figure} [htp]
\centering
\includegraphics[width=1\columnwidth]{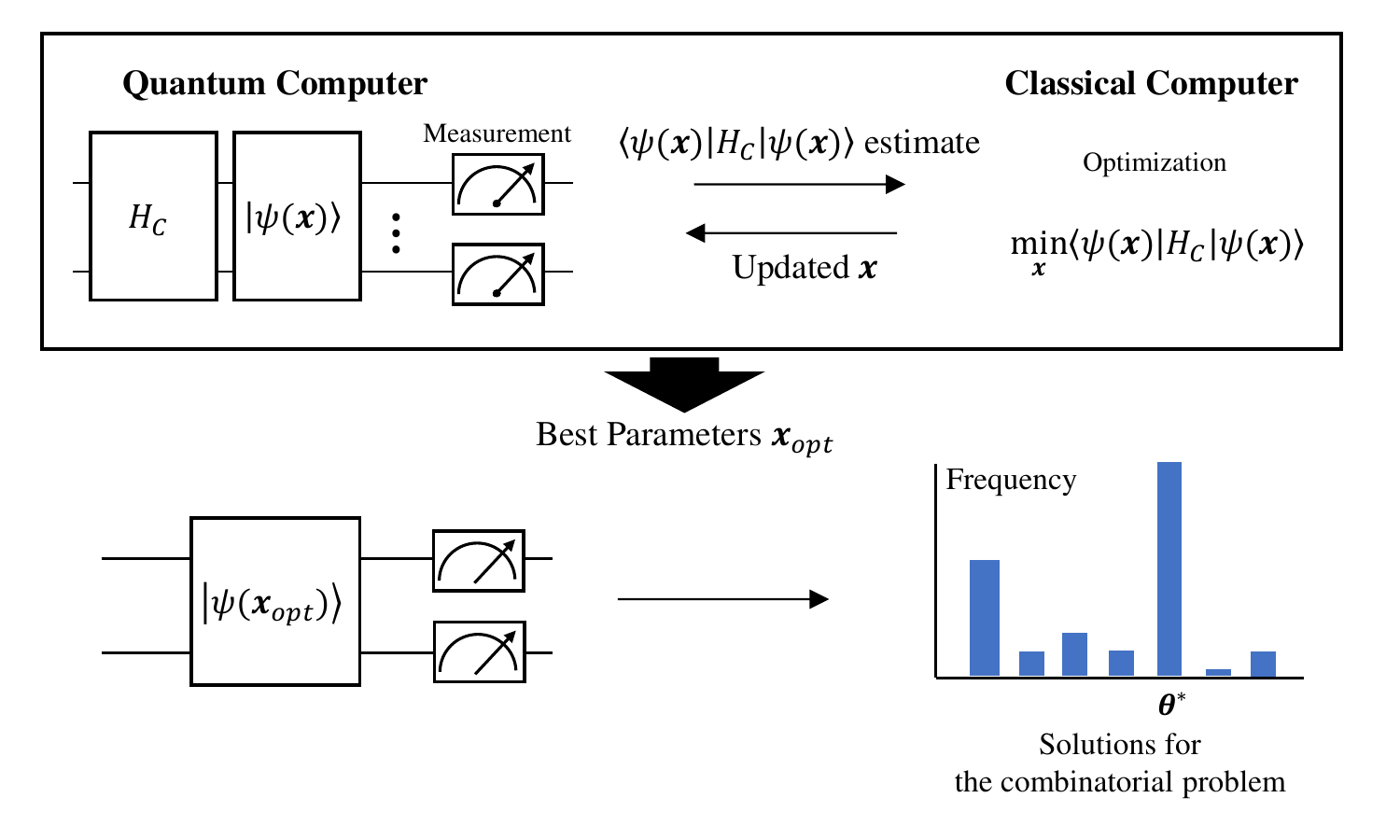}
\caption{In the context of a combinatorial optimization problem with cost Hamiltonian $H_C$, QAOA iteratively updates a parameter vector $\BFx$ to minimize the objective function value in \eqref{eq:qaoa-problem}. Once a sufficiently good solution, $\BFx_{opt}$, is achieved, QAOA proceeds to obtain a probability distribution by measuring the quantum state $\ket{\psi(\BFx_{opt})}$. In this distribution, the solution with the highest frequency corresponds to the optimal solution for the original combinatorial problem, $\bm{\theta}^*$.} \label{fig:qaoa}
\end{figure}

\Cref{fig:qaoa} illustrates the steps of QAOA.
As in our discussion of general VQAs, the quantum computer in \Cref{fig:qaoa} can be viewed as a stochastic oracle that is iteratively queried by a classical computer. 
The parameters $\BFx$ that describe the state $\ket{\psi(\BFx)}$ are updated by the classical computer based on the (stochastic) measurements $\bra{\psi(\BFx)}H_C\ket{\psi(\BFx)}$ made by the quantum computer. 
After some budget is exhausted, or some other stopping criteria  determined by the stochastic optimization algorithm implemented in the classical computer is reached, 
the best observed parameters $\BFx_{opt}$ are measured one more time with a number of shots to yield an empirical discrete distribution on (the finite, but combinatorial, number of) possible bit strings feasible for the combinatorial optimization problem. 
Near the optimal eigenstate solution to \eqref{eq:qaoa-problem}, state vectors have a high probability of collapsing to the optimal solution of the original combinatorial problem upon measurement. 
Therefore, the solution(s) of highest frequency is/are interpreted as candidates for the global optimizer of the combinatorial optimization problem. 

The overall computational expense of executing QAOA can be assessed in a similar way to how SO typically quantifies expense, which is by counting the number of simulation oracle calls (shots) needed to attain a sufficiently accurate solution. 
However, the current state-of-the art for quantum computers involve additional \emph{latencies} that are less seriously considered when designing algorithms for state-of-the-art classical computers. 
Latencies can differ across various architectures. 
For instance, a superconducting quantum processor has measurement times in the range of a few microseconds, see~\cite{gambetta2007protocols}. 
In contrast, a trapped ion system can require hundreds of microseconds to perform a measurement~\citep{clark2021engineering,bruzewicz2019trapped}. 
These measurement times are in addition to the duration of gate operations and system resets, all contributing to the time needed to acquire a single sample of shots.
Moreover, many modern quantum computers operate in a cloud environment, leading to potential extra overhead from network latency, as noted by~\cite{sung2020using}. 
Hence, designing optimization algorithms that consider these latencies is crucial for efficiently using quantum resources in the near term.

The idea of making explicit latency considerations in the design of (theoretical) algorithms for VQAs was made in previous work \citep{menickelly2022quantum}.
The work that we present in this paper is meant to provide a slightly more heuristic, but practical, means to controlling latency within an adaptive sampling framework; this will be seen in our two-stage estimation approach. 

\begin{figure}
    \centering
    \includegraphics[width=.7\textwidth]{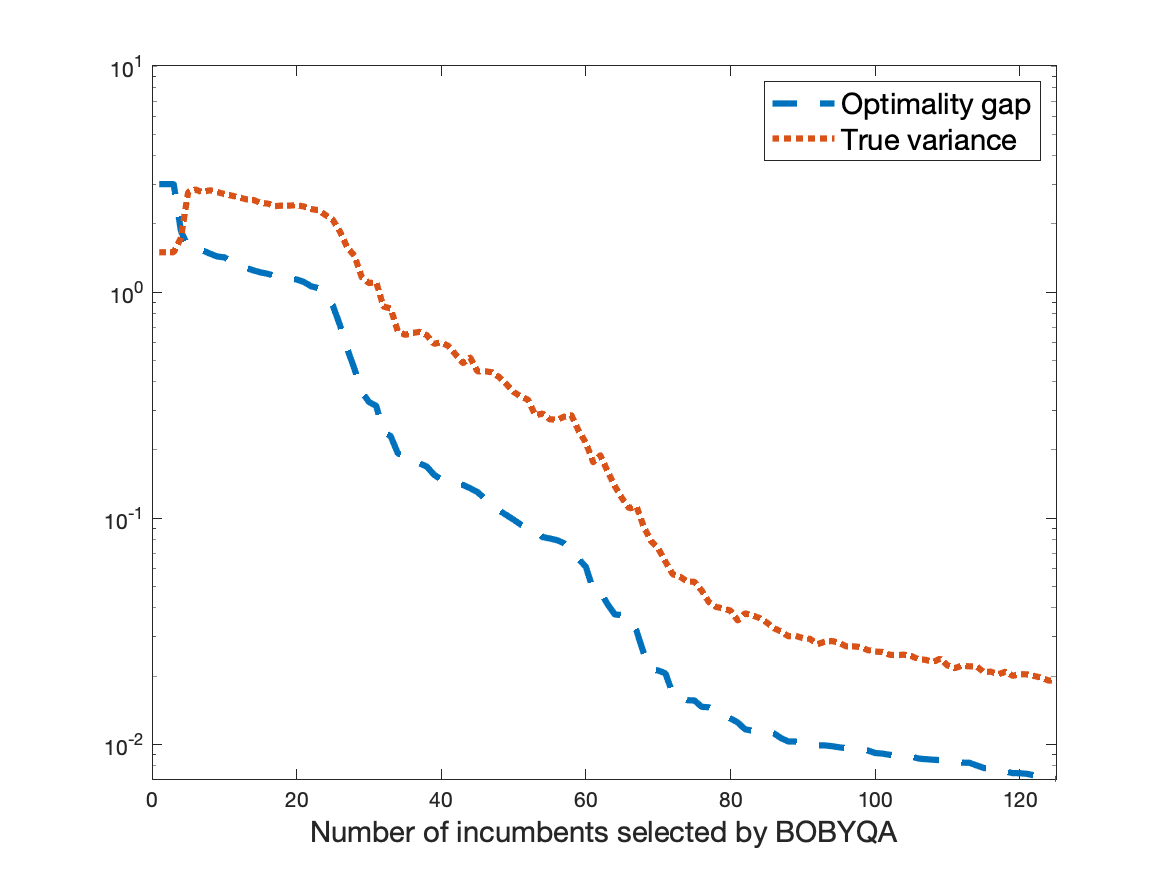}
    \caption{In this plot, we simulate in Qiskit \citep{Qiskit} a depth-10 QAOA circuit for solving a maxcut problem on a toy graph on 6 nodes, for which the optimal solution to the maxcut problem is 6.
    We provide \texttt{BOBYQA} \citep{powell2009bobyqa} with the deterministic expectation statevector value for this toy problem and record the improving sequence of incumbent solutions returned. 
    We illustrate, on the same log scale, the optimality gap of the incumbents found, as well as the population variance associated with the statevector value. 
    We observe that, as expected, population variance decays alongside the optimality gap, but not necessarily monotonically. 
    }
    \label{fig:variance_plot}
\end{figure}

Another distinction of the general VQA setting is that in problem of the form \eqref{eq:qaoa-problem}, it is well-known that an eigenstate should exhibit zero variance; see numerous references within \cite{zhang2022variational}. 
In this paper, we specifically define this characteristic as \emph{state-dependent noise}, represented by the following equation:
\begin{equation}\label{eq:statedependent} \sigma^2(\BFx) \leq C_0 + C_1 (f(\BFx) - f_{\min}), 
\end{equation}
for some $C_0,C_1 \geq 0$, where $\sigma^2(\BFx):=\mbE_{\xi}[(F(\BFx,\xi) - f(\BFx))^2]$ is the true variance of the stochastic function value at $\BFx$ and $f_{\min}$ represents the optimal objective function value. 
While we do not practically require \emph{linearity} in the optimality gap $f(\BFx) - f_{\min}$, as written in \cref{eq:statedependent}, we do coarsely imagine $\sigma^2(\BFx)$ being bounded by \emph{some} function of the optimality gap.  
In~\Cref{fig:variance_plot}, we demonstrate this phenomenon on a small-scale example of using a QAOA circuit for solving a max-cut problem on a toy graph. 
This zero-variance principle has been exploited recently in a quantum computer to self-verify whether a ground state for a given Hamiltonian was accurately prepared \citep{kokail2019self}.
Recently, and also inspired by this phenomenon, \cite{zhang2022variational} considered regularizing VQA cost functions with a measure of estimated variance. 

We further remark that in the VQA setting, there are no  common random numbers that would allow for the faithful reproduction of a sample path, a feature that precludes the use of some techniques in SO.

\subsection{Adaptive Sampling}
In derivative-free SO, stochastic trust-region optimization (TRO) has become increasingly popular as a methodology for solving nonlinear and nonconvex optimization problems~\citep{sun2023tr,cao2022trhighp,jin2021probability,chen2018storm,Sara2018ASTRO,Chang2013STRONG}. 
Stochastic TRO methods generate a random sequence of incumbents, denoted  $\{\BFX_k\}$, during a single run. 
Incumbent selection depends on approximations of the objective function by means of local models, and respective approximate minimizers of these models within dynamically-sized neighborhoods.
In the derivative-free setting, where derivative information is assumed unavailable, these local models are typically computed via interpolation or regression techniques, employing function value estimates at design points near the current incumbent. 
To ensure the accuracy of local models of the objective function, it is imperative to have access to sufficiently accurate function values estimates;
in the stochastic optimization setting, sufficient accuracy can be achieved by averaging a sufficiently large number of samples.
Thus, if an excessive number of simulation oracle calls (shots on a quantum computer) are required to attain this necessary precision, it becomes challenging to find a satisfactory approximate solution to the optimization problem within a reasonable timeframe. 
Hence, for the purpose of judiciously determining a suitable sample size, \cite{Sara2018ASTRO} introduced an adaptive sampling approach within the TRO framework, ASTRO-DF. 
An adaptive sampling strategy dynamically determines the sample size by balancing the estimation error at each point with a measure of first-order optimality error. This dynamic strategy produces a random sample size that is a stopping time with respect to the generated observations at the design point of interest.
 Adaptive sampling has been shown to achieve the highest efficiency measured by the expected total sample complexity~\citep{ha2023,jin2023sample}.
 
 However, and as we set out to address in this paper, the adaptive sampling strategy requires repeated message passing between the optimization engine and the computer simulation, which may be prohibitive in latency-constrained settings. 
One of the key distinctions between the VQA setting and classical SO lies in the fact that in VQA, the ``simulation" is entirely handled by the quantum computer. 
In the VQA setting, a quantum circuit is calibrated according to parameters $\BFx$, which is then executed a number of times. 
Consequently, this process involves a nontrivial amount of communication between the quantum computer and the classical computer, adding an equally nontrivial computational burden to the overall optimization procedure. 
Thus, while an adaptive sampling strategy can reduce the total number of replications by incrementally adding shots until the estimated variance satisfies a particular condition, it may not necessarily alleviate the overall computational burden due to the communication costs. 

We remark that because ASTRO-DF requires an estimate of $\sigma^2(\BFx)$ to compute the number of samples (shots) requested at a design point $\BFx$, and because this number of shots scales linearly with the estimate of $\sigma^2(\BFx)$, we anticipate that in state-dependent noise settings, near optimality, this decaying variance will play a mitigating role in the number of required samples, as contrasted with the decreasing trust-region radius. 

\subsection{Our Contributions} 
Motivated by the particularities of VQA problems, 
we propose 1) a replacement for the adaptive sampling strategy in ASTRO-DF with a two-stage sampling strategy, and 2) a refinement for the local model construction employed in ASTRO-DF. 
Both of these contributions hinge on 
 a secondary model that interpolates or regresses variance estimates of the stochastic objective function, in order to locally approximate the variance of previously unevaluated design points or incumbents. 
The proposed two-stage estimation approach will ensure at most two communications between the optimization engine and the quantum computer per function evaluation. The variance model helps to achieve this by predicting the variance at previously unevaluated parameters. 

Predictions from the variance model will additionally aid in choosing interpolation points in new incumbent neighborhoods. 
In particular, when selecting design points for objective model-building, we will prioritize points exhibiting lower model variance.  
A particular advantage of this prioritization is that it provides a heuristic intended to escape local minima in objective functions exhibiting state-dependent noise.

We will delve into the details of our proposed uses of a variance model throughout the optimization process in \Cref{sec:point-selection} and \Cref{sec:how-ss}, respectively.
Taken together, these improvements are intended to alleviate the insistence, in the usual adaptive sampling setting, on incremental sampling to estimate if the estimated variance at a previously unevaluated design point or incumbent is sufficiently small. 
Instead, we begin the adaptive sampling step with a reasonable initial estimate of estimated variance gleaned from the auxiliary model, and then only require at most two additional communications with the quantum device to further refine the estimate. 

\section{Simulation Optimization with Trust Regions}
Stochastic TRO is effective at solving zeroth-order nonconvex stochastic optimization problems. Its salient feature is a natural ability to self-tune step sizes and facility for incorporating approximate curvature information.
We provide a fairly generic framework that describes the vast majority of derivative-free stochastic TRO methods. 
A set of design points $\mcX_k$ are evaluated in a neighborhood of the incumbent $\BFX_k$.
A model $M_k(\cdot)$ is fit to those function evaluations at $\mcX_k$.
An approximate minimizer of $\mcX_k$ over a trust region of size $\Delta_k$, i.e., $\mcB(\BFX_k;\Delta_k)$, is computed and denoted $\BFXtilde_{k+1}$. 
If $\BFXtilde_{k+1}$ witnesses sufficient decrease over $\BFX_k$, then we set the next incumbent as $\BFX_{k+1} = \BFX_k$.
ASTRO-DF is a variant of this class of algorithms that embeds adaptive sampling to determine a judicious lower bound on the number of oracle calls (shots) required at each design point to guarantee optimization progress. 
The key element of this approach involves allocating computational resources based on a measure of the optimality gap, such as $\|\nabla f(\BFX_k)\|$, which ASTRO-DF consistently monitors by means of the trust-region radius $\Delta_k$. 
As a result, it is typical that more computational effort is expended on points in closer proximity to first-order critical regions. 
Before we expound on recent developments in ASTRO-DF, we begin by introducing the notation and definitions that will be employed throughout this paper.

\subsection{Notation and Definition}\label{sec:notation}
We will use capital letters for random variables, bold fonts for vectors, script fonts for sets and $\sigma$-fields, and a sans-serif font for matrices.

In addition to $f(\BFx)$ and $F(\BFx,\xi)$, we additionally define a sample average based on $N(\BFx)$ many samples
$$\bar{F}(\BFx, N(\BFx)) = \frac{1}{N(\BFx)} \displaystyle\sum_{i=1,\dots,N(\BFx)} F(\BFx, \xi_i),$$
where $\{\xi_i: i=1,2,\dots, N(\BFx)\}$ denote independent realizations of $\xi_i$. 

\begin{definition}(stochastic interpolation models).
    Let $\Phi(\BFx)=(\phi_0(\BFx),\phi_1(\BFx), \dots, \phi_q(\BFx))$ form a linearly independent set of polynomials on $\mbR^d$. 
    With $q = p+1$, $\BFX_k^{0}:=\BFX_k$, and the design set $\mcX_k:=\{\BFX^{i}\}_{i=1}^{p}\subset \mcB(\BFX_k;\Delta_k)$, consider the linear system
    \begin{equation}\label{eq:syslineq}
        \sfM(\Phi, \mcX_k) \BFbeta_k = \Fbar(\mcX_k,N(\mcX_k)),
    \end{equation}
    where \\ 
    \begin{equation*}
        \sfM(\Phi, \mcX_k) = 
        \begin{bmatrix}
        \phi_1(\BFX_k^{0}) & \phi_2(\BFX_k^{0}) & \cdots & \phi_q(\BFX_k^{0})  \\
        \phi_1(\BFX_k^{1}) & \phi_2(\BFX_k^{1}) & \cdots & \phi_q(\BFX_k^{1})  \\
        \vdots & \vdots & \vdots & \vdots \\
        \phi_1(\BFX_k^{p}) & \phi_2(\BFX_k^{p}) & \cdots & \phi_q(\BFX_k^{p})  \\
        \end{bmatrix},   
        \Fbar(\mcX_k,N(\mcX_k)) =
        \begin{bmatrix}
        \Fbar(\BFX_k^{0},N(\BFX_k^{0})) \\
        \Fbar(\BFX_k^{1},N(\BFX_k^{1})) \\
        \vdots \\
        \Fbar(\BFX_k^{p},N(\BFX_k^{p})) \\
        \end{bmatrix}.
    \end{equation*}

    We say the set $\mcX_k$ is \emph{poised with respect to $\Phi(\BFx)$} if the matrix $\sfM(\Phi, \mcX_k)$ is nonsingular. 
    If there exists a solution $\BFbeta_k = (\beta_{k,i},\ i=0,1,2,\ldots,p)$ to~\eqref{eq:syslineq} (i.e., $\mcX_k$ is poised with respect to $\Phi(\BFx)$),   
    then the function $M_k:\mcB(\BFX_k;\Delta_k) \to \mbR$, defined as $M_k(\BFx) = \sum_{j=0}^{p} \beta_{k,j} \phi_{j}(\BFx)$ is a \emph{stochastic interpolation model} of estimated values of $f$ on $\mcB(\BFX_k;\Delta_k)$. 
\label{defn:polyintermd}
\end{definition}

We specialize~\Cref{defn:polyintermd} to the case where a diagonal model Hessian is employed. 

\begin{definition} (stochastic quadratic interpolation models with diagonal Hessian)
Let $\Phi(\BFx)$ be the polynomial basis $\{1, x_1,\dots, x_d, x_1^2,\dots, x_d^2\}$ so that $p = 2d$. 
Let $\mcX_k$ be poised with respect to $\Phi(\BFx)$ and let $M_k(\BFx)$ be a stochastic interpolation model.
Then, 
denoting $\BFG_k:=\begin{bmatrix}
    \beta_{k,1} & \beta_{k,2} & \cdots & \beta_{k,d} 
    \end{bmatrix}^\intercal$ 
and letting 
$\sfH_k$ be a $d\times d$ matrix with $[\sfH_k]_{i,i} = \beta_{k,d+i}$ and zeros off the diagonal, 
we refer to 
\begin{equation}
    M_k(\BFx) = \beta_{k,0} +  (\BFx-\BFX_k)^\intercal \BFG_k + \frac{1}{2} (\BFx-\BFX_k)^\intercal \sfH_k(\BFx-\BFX_k),\label{eq:mdefn}
\end{equation}
as a 
\emph{stochastic quadratic interpolation model of $f$ with a diagonal Hessian}. 
\label{defn:diagHess}
\end{definition}
 
A particular utility of \Cref{defn:diagHess} is that any coordinate stencil, such as 
$\mcX_k^{cb} = \{\BFX_k, \BFX_k+\BFe_1 \Delta_k, \ldots, \BFX_k+\BFe_d \Delta_k, \BFX_k-\BFe_1 \Delta_k, \ldots, \BFX_k-\BFe_d \Delta_k \}$ where 
$\BFe_i$ denotes the $i$-th elementary basis vector of $\mbR^d$, 
is clearly poised with respect to $\Phi(\BFx)$; \Cref{defn:diagHess} is therefore immediately nonvacuous. 

\begin{definition} (filtration) A \emph{filtration $\{\mcF_{k}\}_{k \geq 1}$ over a probability space $(\Omega,\mbP,\mcF)$} is an increasing sequence of $\sigma$-algebras within $\mcF$, where each $\mcF_k$ is a subset of $\mcF_{k+1}$, and all are subsets of $\mcF$, for every $k$. $\mcF_k$ represents all information that is available at time $k$. 
\end{definition}

\subsection{History-informed ASTRO-DF}
Recent augmentations to ASTRO-DF have aimed at boosting computational efficiency. We collectively refer to ASTRO-DF with these augmentations as ``history-informed ASTRO-DF". History-informed ASTRO-DF includes a direct search (see, e.g., \cite{ha2023jsim,ha2023wsc}) component in each iteration to increase the likelihood of finding a new incumbent in each iteration without increasing the allotted budget. 
In cases where $\BFXtilde_{k+1}$, the local model minimizer, does not lead to a sufficient reduction in the estimated function value, ASTRO-DF would have declared the iteration as unsuccessful, immediately contracting the trust-region radius.
This declaration of an unsuccessful iteration and trust-region contraction would have occurred even if some design points in $\mcX_k$ yielded improvements over the incumbent. 
Electing to replace the next iterate with the best design point from $\mcX_k$ is tantamount to a direct search iteration. 
Practically, having fewer unsuccessful iterations due to the direct search feature amounts to a slower rate of decay in the trust-region radius $\Delta_k$; in turn, this slower rate keeps the sample size (which is proportional to $\Delta_k^{-4}$--see \eqref{eq:sample_condition}) from growing too quickly. 

Moreover, history-informed ASTRO-DF consistently employs a (rotated) coordinate stencil like $\mcX_k^{cb}$ as the set of design points via a reuse strategy to conserve computational resources in each iteration. 
If there are previously evaluated design points located within the trust region in an iteration, 
the design point that is farthest from the incumbent $\BFX_k$ (denote it $\BFY$) will be added to the design set $\mcX_k^{cb}$, and the simulation results at that point will be reused. 
The remaining members of the design set $\mcX_k^{cb}$ will be selected deterministically by computing a set of mutually orthonormal vectors, all orthogonal to $\BFY - \BFX_k$ and using these vectors to generate a rotated coordinate basis.
The benefit of employing a rotated coordinate basis is that we obtain a more precise gradient estimate at $\BFX_k$ than if we had used alternative bases \citep{ha2023wsc}; in fact, among all design sets satisfying $|\mcX_k| = 2d+1$, a rotated coordinate basis is optimal in a precise sense \citep{tom2023optimalpoised}.
As a practical matter, it is often the case that $\BFY = \BFX_{k-1}$ following a successful iteration; therefore, this reuse strategy, in tandem with the direct search strategy, will extrapolate the search direction from the previous iteration. 
\Cref{tab:astrodf} summarizes several key differences between the original ASTRO-DF and history-informed ASTRO-DF. 

\begin{table}[!htbp]\caption{{Differences between ASTRO-DF and history-informed ASTRO-DF.}} 
    \label{tab:astrodf}
    \centering
    \begin{tabular}{|c | c | c | c|} 
        \hline
        Algorithm & ASTRO-DF & History-informed ASTRO-DF   \\ [0.5ex] 
        \hline\hline
        Selection of $\mcX_k$ & Random & Rotated coordinate basis  \\ \hline
        $|\mcX_k|$ & $(d+1)(d+2)/2$ & $2d+1$ \\ \hline
        Source of next incumbent & Model & Model and $\mcX_k$  \\ \hline
        $|\mcX_k\cap \mcX_{k-1}|$ & $\ge 0$ &  2 \\ \hline
    \end{tabular}
\end{table}


\section{Point Selection}
\label{sec:point-selection}
The efficiency of a TRO algorithm is closely tied to the geometry of the design set $\mcX_k$ in each iteration, since the choice of $\mcX_k$ directly impacts the quality of the local model. 
As an abstract example, if the design points in $\mcX_k$ lie entirely in a single halfspace intersected with the trust region, the local model $M_k$ is prone to having poor predictive accuracy on the complementary halfspace. 
In addition to geometry, the quantity of design points contained in $\mcX_k$ in each iteration plays a crucial role. 
Employing an excessive number of design points per iteration can lead to prolonged computation times. 
Conversely, if we use too few design points relative to the choice of basis $\Phi(\BFx)$ in \Cref{defn:polyintermd} so that $p + 1 < q$, then the system \eqref{eq:syslineq} is underdetermined, and local models satisfy error bounds for underdetermined systems, that are generally worse than the error bounds derivable for determined systems of equations \eqref{eq:syslineq}; see e.g., \citep{katya:DFObook}[Section 5]. 
To address all of these issues, and based on practical experience, history-informed ASTRO-DF makes the explicit choice to select $2d+1$ design points per iteration and always employs a rotated coordinate basis, see \Cref{defn:diagHess}.
See \Cref{fig:rcb} for an illustration of the rotated coordinate basis selection performed by history-informed ASTRO-DF. 
This choice of basis and determined $\mcX_k^{cb}$ allows us to capture some curvature information within a reasonable computational timeframe, and simultaneously realizes a well-distributed placement of design points within the trust region to faithfully approximate the objective function. 
Past empirical evidence has demonstrated that history-informed ASTRO-DF is superior in performance to the older implementation of ASTRO-DF that makes a more general choice of basis \citep{ha2023wsc} on a range of SO problems, as measured by progress made within a simulation budget. 

\begin{figure} [htp]
\centering
\subfloat[Original point selection strategy using $M_k$]{%
\resizebox*{8cm}{!}{\includegraphics{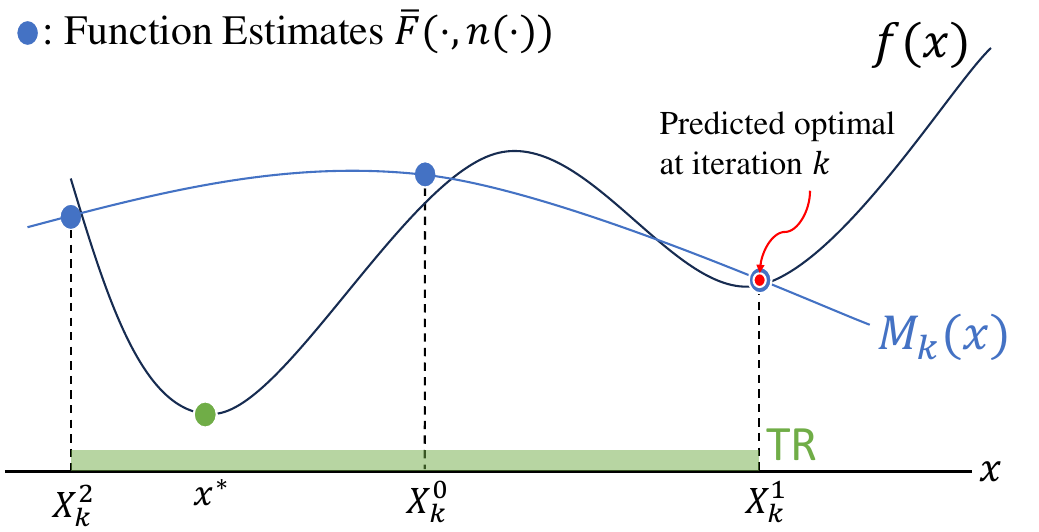}}\label{fig:wovm}}
\subfloat[New point selection strategy using both $M_k$ and $M^v_k$.]{%
\resizebox*{8cm}{!}{\includegraphics{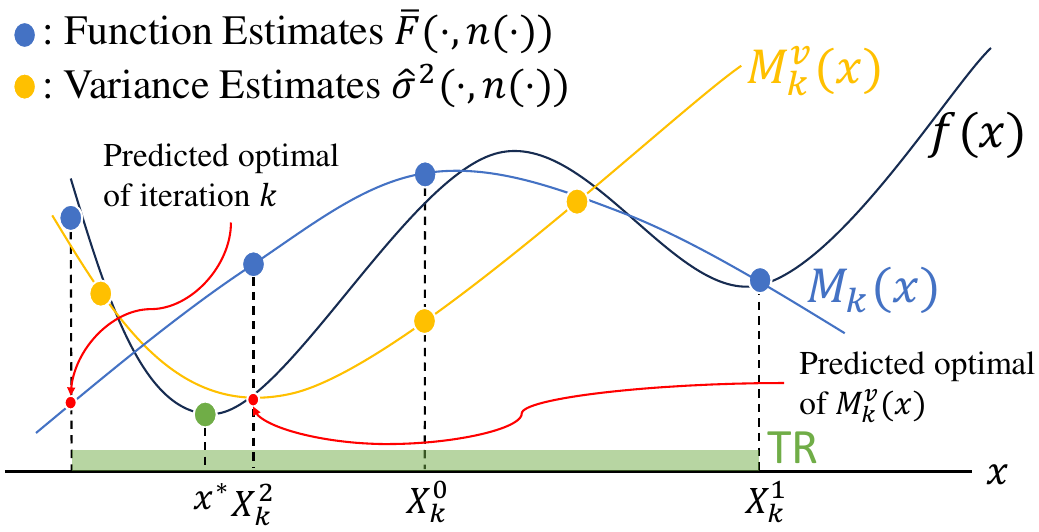}}\label{fig:wvm}}
\caption{A cartoon illustrating the effect of using the two local models. \Cref{fig:wovm} illustrates the performance of history-informed ASTRO-DF without the inclusion of $M_k^v$. In this case, $\BFXtilde_{k+1}$ will be further drawn to the basin of a local minimum. \Cref{fig:wvm} illustrates the performance of history-informed ASTRO-DF with the proposed use of two models, $M_k$ and $M_k^v$. In this case, the design point $\BFX_k^2$, the minimizer of $M_k^v$, will be selected as the next incumbent $\BFX_{k+1}$ We prefer this global optimality-seeking behavior. }
\label{fig:with-var-model} 
\end{figure}

In this paper, we propose exploiting a previously discussed characteristic of VQAs, namely, the diminishing variance with respect to the optimality gap. 
To harness the potential acceleration of ASTRO-DF through the utilization of this distinctive property, we propose an extension of the point selection strategy of history-informed ASTRO-DF.
Central to the new point selection strategy is the construction of a local model that interpolates or regresses variance estimates. 
Within ASTRO-DF, upon computation of sample average function values at each design point in $\mcX_k$, we also acquire estimated variance, i.e., $\widehat{\text{Var}}(F(\cdot,\xi))$, at each design point in $\mcX_k$ with little additional computational overhead; this overhead is arithmetic, and not incurred by the simulation oracle (quantum computer). 
Then, without additional oracle effort, we can construct a local model of variance by solving \eqref{eq:syslineq}, where we replace the right hand side with estimated variance, as opposed to estimated function value. We denote this model by $M_k^v(\BFx)$. 
We denote the (approximate) minimizer of $M_k^v(\BFx)$ within the trust region by $\BFXtilde^v_k$.
Coupled with the discussed fact that history-informed ASTRO-DF utilizes a direct search in each iteration, if the design set $\mcX_k$ includes $\BFXtilde^v_k$, 
then we will end up considering $\BFXtilde^v_k$ as a potential next incumbent $\BFX_{k+1}$. 
Owing to the correlation between variance and optimality gap frequently observed in VQAs, 
this proposed scheme provides a heuristic intended to drive the incumbents $\BFX_k$ to \emph{global} optimality, since the global minimizer for $f$ should coincide with the minimizer of the true variance.
We stress that this proposed scheme is indeed a heuristic, since its utility is also a function of the trust-region radius, but it is a seemingly useful heuristic. 
\Cref{fig:with-var-model} shows an example demonstrating the effect of using two local models $M_k$ and $M_k^v$.
We summarize what we refer to as the \emph{two-model approach} in the following algorithm and illustrate one iteration of the two-model approach in \Cref{fig:wvm}:
\begin{enumerate}
    \item Construct the local model $M_k^v(\BFx)$ using previously evaluated design points and their associated sample variance estimates and find the minimizer, $\BFXtilde^v_k$ of $M_k^v(\BFx)$ within the trust region.
    \item Determine the design set $\mcX_k$, containing $\BFX^v_k$, and estimate the function at each point in $\mcX_k$.
    \item Construct the local model $M_k$ using $\mcX_k$ and the corresponding function value estimates, then minimize $M_k$ within the current trust region to obtain $\BFXtilde_{k+1}$. 
    \item Obtain a function value estimate at $\BFXtilde_{k+1}$ and then utilize the direct search on all points in $\mcX_k \cup \{\BFXtilde_{k+1}\}$ to determine the next incumbent $\BFX_{k+1}$. 
\end{enumerate}

\begin{figure} [htp]
    \centering
    \includegraphics[width=.7\textwidth]{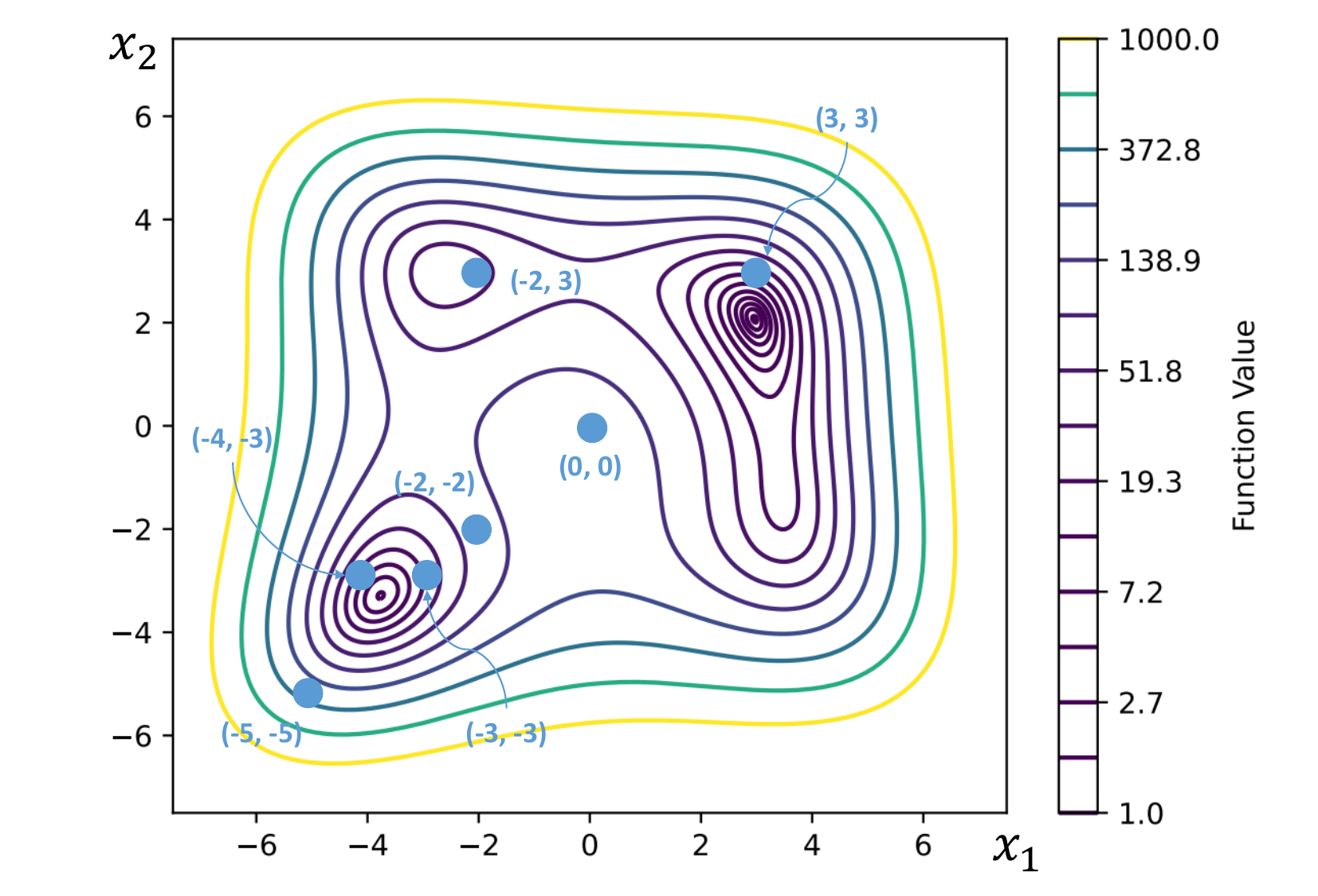}
    \caption{Contour plot of the expectation of \eqref{eq:himme-problem}. Each labelled point corresponds to an experimented initial solution shown in Figures~\ref{fig:himme-diff-init} and \ref{fig:himme-diff-(3,3)}. The global minimum is attained at $(3, 2)$.
    }
    \label{fig:contour}
\end{figure}

We now will demonstrate the efficacy of this two-model approach by performing some preliminary experiments designed to highlight the utility of seeking variance-minimizing design points. 
The test problem involves a stochastic variant of the Himmelblau function.
We intentionally devise the stochastic noise so as to be state-dependent; 
in particular, we make the variance vanish at the global minimum of the Himmelblau function.
The test problem is
\begin{equation} \label{eq:himme-problem}
F(\BFx,\xi) = (x_1^2 + x_2 - 11)^2 + (x_1 + x_2^2 -7)^2 + |x_1-3| +\xi,    
\end{equation}
where $\xi~\mcN(0,|(x_1-3)(x_2-2)|)$. 
The global minimum for the noiseless version of the Himmelblau function is located at $(3, 2)$; by construction, the variance is 0 at $(3,2)$, while the other local optima exhibit positive.
The expectation of \eqref{eq:himme-problem} is depicted in \Cref{fig:contour}.

\begin{figure} [htp]
\centering
\subfloat[initial point = (-5,-5)]{%
\resizebox*{5.2cm}{!}{\includegraphics{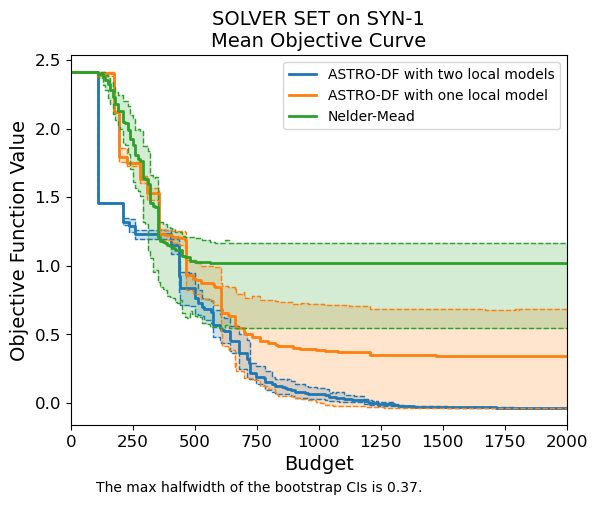}}\label{fig:(-5,-5)}}
\subfloat[initial point = (0,0)]{%
\resizebox*{5.2cm}{!}{\includegraphics{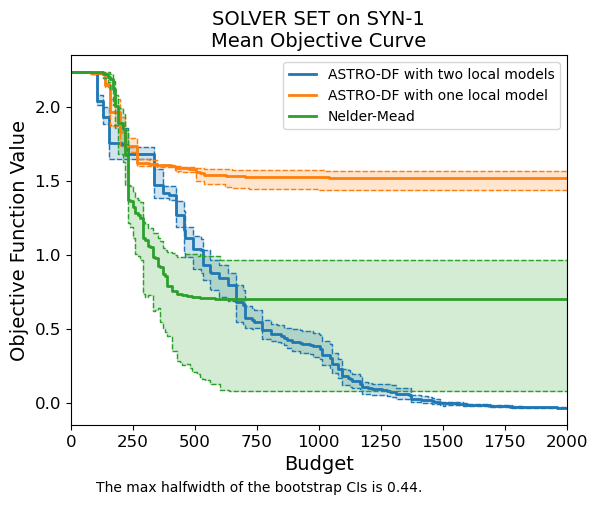}}\label{fig:(0,0)}} 
\subfloat[initial point = (-2,-2)]{%
\resizebox*{5.2cm}{!}{\includegraphics{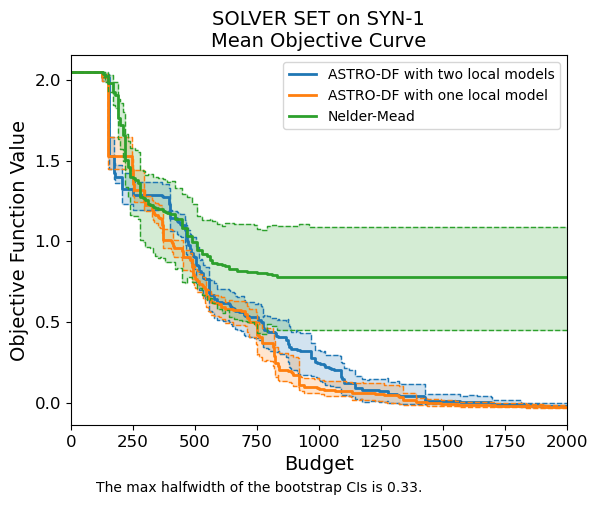}}\label{fig:(-2,-2)}}
\\
\subfloat[initial point = (-4,-3)]{%
\resizebox*{5.2cm}{!}{\includegraphics{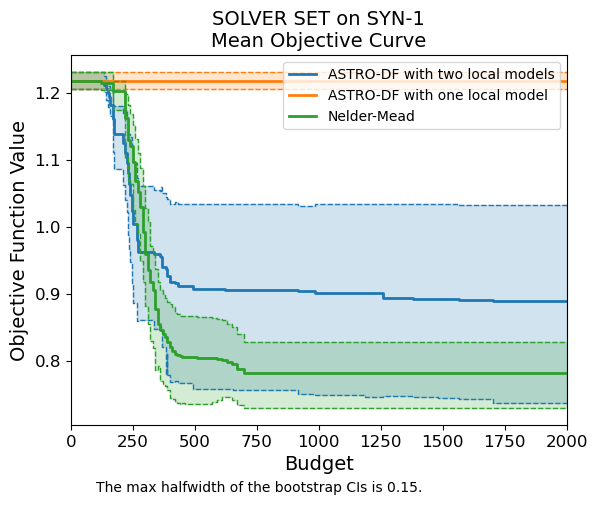}}\label{fig:(-4,-3)}} 
\subfloat[initial point = (-3,-3)]{%
\resizebox*{5.2cm}{!}{\includegraphics{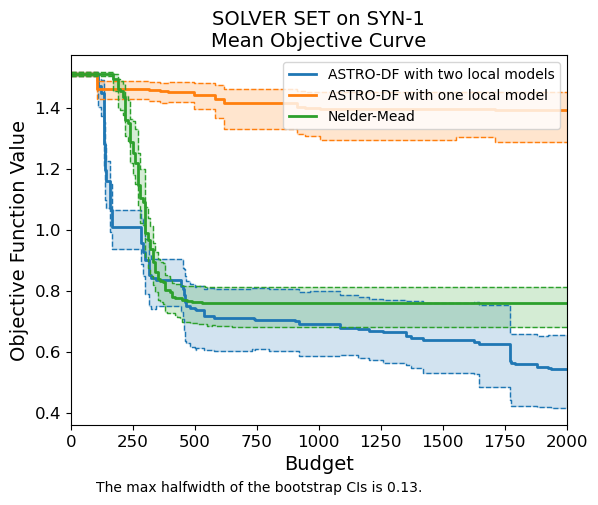}}\label{fig:(-3,-3)}}
\subfloat[initial point = (-2,3)]{%
\resizebox*{5.2cm}{!}{\includegraphics{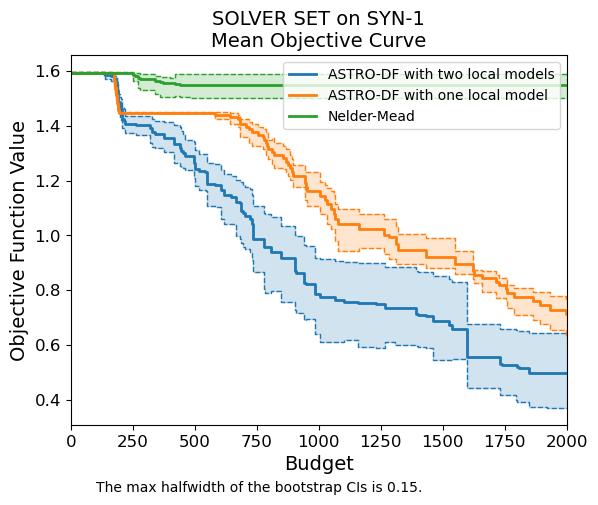}}\label{fig:(-2,3)}}\\

\caption{Performance of solvers on \eqref{eq:himme-problem} provided with various initial points. 
Translucent bands represent a 95\% confidence interval over 20 macro-replications and solid lines represent mean performance. provided with various initial points. The $y$-axis is on a logarithmic scale.}
\label{fig:himme-diff-init}
\end{figure}

We begin our investigation by examining whether two-model ASTRO-DF can improve the likelihood of discovering the global optimum when compared directly to what we will refer to as \emph{single-model} ASTRO-DF, that is, the history-informed ASTRO-DF. 
In \Cref{fig:himme-diff-init}, we initialize one-model and two-model ASTRO-DF (and Nelder-Mead as a baseline comparison) from various initial solutions and plot the best (estimated) objective function values attained as a function of oracle calls. To ensure robustness and accuracy, each algorithm was executed 20 times, referred to as 20 macro-replications. 
From each initial solution, two-model ASTRO-DF exhibits a higher probability of identifying the basin of $(3,2)$ than single-model ASTRO-DF.
Notably, when the initial design point is relatively far from any local minimum, it becomes evident that two-model ASTRO-DF excels at discovering the global minimum, as evident in the observations from Figure \ref{fig:(-5,-5)} and \ref{fig:(0,0)}. 
Furthermore, it is noteworthy that while some runs may be slow to converge to the global minimum when the optimization process is initialized near a local optimum, the majority of them still manage to reach the global minimum eventually, as illustrated in Figure \ref{fig:(-4,-3)}, \ref{fig:(-3,-3)}, and  \ref{fig:(-2,3)}.
Finally, when the initial design point is placed relatively near the global minimum, the convergence rate can be relatively slow, see Figure~\ref{fig:(3,3)} and \ref{fig:(3,3)-50000}. 
This is because the trust region must become sufficiently small to identify a new incumbent; nonetheless, with a sufficiently large budget (Figure~\ref{fig:(3,3)-50000}), ASTRO-DF eventually does find decrease. (Sensitivity to initial trust region size is a general TRO criticism, outside our scope here.)

\begin{figure} [htp]
\centering
\subfloat[small budget]{%
\resizebox*{6cm}{!}{\includegraphics{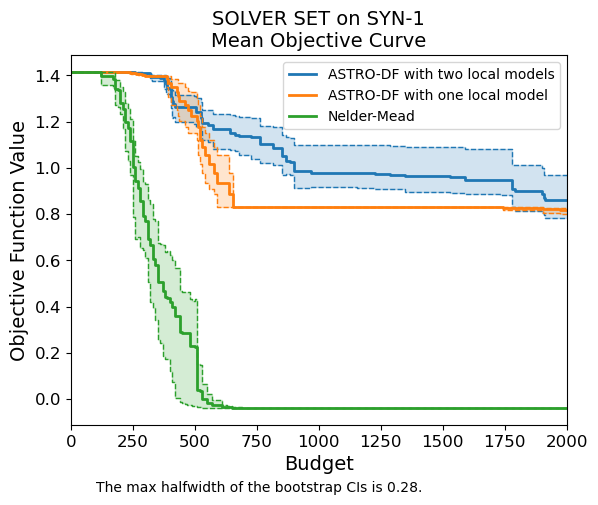}}\label{fig:(3,3)}}\hspace{12pt}
\subfloat[large budget]{%
\resizebox*{6cm}{!}{\includegraphics{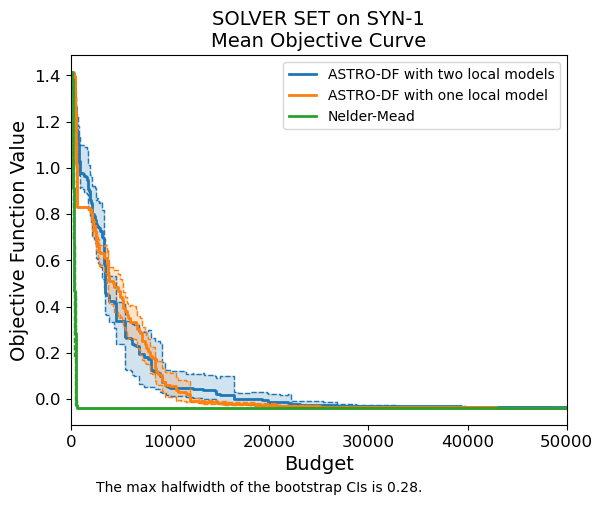}}\label{fig:(3,3)-50000}} 
\caption{Performance on \eqref{eq:himme-problem} with 95\% confidence interval. The $y$-axis is again on a logarithmic scale. The initial point, $(3, 3)$, is relatively close to the global minimum.}
\label{fig:himme-diff-(3,3)}
\end{figure}

\section{Sample Size Selection}
\label{sec:how-ss}
In the original implementation of ASTRO-DF \citep{Sara2018ASTRO}, the sample size at any design point $\BFx \in \mbR^d$ is determined by the formula
\begin{equation}
\label{eq:sample_condition}
    N_k(\BFx)=\min\biggl\{ n\geq \lambda_k:\frac{\sigmahat\left(\BFx,n\right)}{\sqrt{n}}\leq\frac{\kappa\Delta_{k}^{2}}{\sqrt{\lambda_k}}\biggr\} \text{ for } i=0,1,\ldots,2d,
\end{equation}  
where $\{\lambda_k\}$ represents a deterministically increasing sequence with logarithmic growth that represents the minimum sample size at iteration $k$, $\sigmahat^2\left(\BFx,n\right)$ denotes the estimated variance from $n$ samples (observations), and $\kappa>0$ is a user-defined constant. 
The right-hand side of the inequality in 
\eqref{eq:sample_condition} is a mildly deflated proxy for the true gradient norm 
(a measure of the first-order optimality gap).
Enforcing a lock-step between estimation error and optimization error ($\kappa\Delta_k^2$) in \cref{eq:sample_condition} helps identify a ``just right" (neither too small nor too large) sample size, and hence enhances algorithm efficiency. 
This ideal sample size (provided it is larger than $\lambda_k$) then satisfies
\begin{equation}
    N_k(\BFx)\geq \lambda_k\frac{\sigma^2(\BFx)}{\kappa\Delta_k^4}.\label{eq:as}
\end{equation} 
However, because $\sigma^2(\BFx)$ is unknown, its most recent estimate $\sigmahat^2(\BFx,n)$ can inform whether the estimated function value with sample size $n$ satisfies \eqref{eq:as}.
In the original implementation of ASTRO-DF, the sample size $n$ is increased by 1 (or a small batchsize) until \eqref{eq:sample_condition} is satisfied. 
In other words, $N_k(\BFx)$ is a stopping time random variable that is learned on the fly, adapting tightly to new observations.
At stopping, there is quantifiable certainty that $N_k(\BFx)$ is the first sample size satisfying
\begin{equation}
    N_k(\BFx)\geq \lambda_k\frac{\sigmahat^2(\BFx,N_k(\BFx))}{\kappa\Delta_k^4}.\label{eq:as2}
\end{equation}

As discussed, when faced with the issue of latency in accessing a quantum device, as in the VQA setting, this ``streaming" approach to adaptive sampling becomes costly, since each pass through the streaming loop will incur additional communication costs.  
To address this issue, we propose three different two-stage estimation strategies. 
The fundamental idea behind each of these three proposed strategies centers around restricting the number of accesses to the quantum computer to at most two. 
The distinction among these three strategies lies in the approach taken to choose a sample size for an initial variance estimate at the point $\BFx$ being evaluated. 

\subsection{Two-stage estimation using initial variance estimate with sample size $\lambda_k$} We first consider simply using $\lambda_k$ as the initial sample size. 
In this case, irrespective of the actual variance at $\BFx$, the first-stage sample size is $N_{k,1}(\BFx) = \lambda_k$; the corresponding estimated variance will be simply $\sigmahat^2(\BFx,\lambda_k)$. 
Then, the second-stage sample size will be computed as 
\begin{equation} \label{eq:ss-lambda}
    N_{k}(\BFx) = \lambda_k\max\left\{1,\dfrac{\sigmahat^2(\BFx,\lambda_k)}{\kappa \Delta_k^4}\right\}, 
\end{equation}
which is derived from~\eqref{eq:sample_condition}. This two-stage process is summarized in \Cref{alg:two-lambda}. 

\begin{algorithm}[htp]  
\caption{$N_k(\BFx)$=\texttt{TwoStageEstimation($\Delta_k,\BFx,\kappa,\mcF_k$)}}
\begin{algorithmic}[1]
\label{alg:two-lambda}
\REQUIRE trust-region radius $\Delta_k$, design point $\BFx$, and minimum sample size $\lambda_k$.
\STATE Set $N_{k,1}(\BFx) = \lambda_k$ and evaluate the estimate $\sigmahat(\BFx,N_{k,1}(\BFx)).$
\STATE Return $N_k(\BFx)$ using \eqref{eq:ss-lambda}
\end{algorithmic}
\end{algorithm}

If $\lambda_k$ is sufficient to attain an accurate function and variance estimate, that is $N_{k}(\BFx) \leq \lambda_k$, then there is no need for additional simulation oracle calls, and the function estimate at $\BFx$ will be $\Fbar(\BFx,\lambda_k)$.
However, with this strategy, when $\lambda_k$ is excessively small, $\sigmahat^2(\BFx,\lambda_k)$ may be a poor estimator of the variance. 
This becomes particularly noticeable when $\sigmahat^2(\BFx,\lambda_k)$ is very large, potentially resulting in an unnecessarily large value for $N_{k}(\BFx).$ 
Given that we are motivated by practicality and limited oracle budgets, this can greatly disrupt our goal of achieving a sufficiently good solution within a reasonable timeframe.
The next strategy is proposed to ameliorate this issue. 

\subsection{Two-stage estimation using variance estimates from $M_k^v$} 
We double-purpose the variance model $M_k^v$ described in \Cref{sec:point-selection} for deriving a two-stage estimation scheme. 
That is, the first-stage sample size $N_{k,1}$ will be determined using the predicted variance at $\BFx$ via the variance model, i.e., 
\begin{equation} \label{eq:initial-ss-var-approx}
    N_{k,1}(\BFx) = \lambda_k\max\left\{1,\dfrac{M_k^v(\BFx)}{\kappa \Delta_k^4}\right\}.
\end{equation}
 We then compute the estimated variance $\sigmahat^2(\BFx,N_{k,1}(\BFx))$ and re-adjust the sample size, if needed, for the second stage in a bid to meet the \eqref{eq:as2} criterion. 
This two-stage estimation technique is summarized in Algorithm~\ref{alg:two-variance-mod}.

\begin{algorithm}[htp]  
\caption{$N_k(\BFx)$=\texttt{TwoStageEstimation($\Delta_k,\BFx,M_k^v,\kappa,\mcF_k$)}}
\begin{algorithmic}[1]
\label{alg:two-variance-mod}
\REQUIRE trust-region radius $\Delta_k$, design point $\BFx$, variance model $M_k^v$, and history $\mcF_k$.
\STATE Set $N_{k,1}(\BFx)$ using \eqref{eq:initial-ss-var-approx} and evaluate the estimate $\sigmahat(\BFx,N_{k,1}(\BFx)).$
\STATE Return 
\[N_{k}(\BFx) = \max\left\{N_{k,1}(\BFx),\lambda_k\dfrac{\sigmahat^2(\BFx,N_{k,1}(\BFx))}{\kappa \Delta_k^4}\right\}.\]
\end{algorithmic}
\end{algorithm}

This strategy presents a new challenge – if the accuracy of the local model $M_k^v$ is poor, $M_k^v(\BFx)$ can be excessively large, again leading to an unnecessary overexpense to our budget.
This consideration motivates our final strategy, which hybridizes Algorithms~\ref{alg:two-lambda} and \ref{alg:two-variance-mod}. 

\subsection{Hybrid two-stage estimation.} 
In this final proposed strategy, our general preference is to utilize $M_k^v$ over $\lambda_k$, as in \Cref{alg:two-variance-mod}. 
However, in cases where $M_k^v$ proves to be inaccurate, we would prefer to pivot to employing the estimated variance with $\lambda_k$ many samples instead, as in \Cref{alg:two-lambda}. 
Thus, we require some heuristic to decide whether or not we trust $M_k^v$ in each iteration. 
To derive this heuristic, we make an assumption that the variance function is Lipschitz continuous within the trust region, i.e., \[|\sigma^2(\BFx_1)-\sigma^2(\BFx_2)|\leq L_v\|\BFx_1-\BFx_2\|,\ \forall \BFx_1,\BFx_2\in\mbR^d,\]where $L_v$ represents a Lipschitz constant. 
This Lipschitz continuity implies that when $M_k^v(\BFx)$ exceeds the variance estimate at the current iterate, $M_k^v(\BFX_k)$, by more than a constant factor times $\|\BFx - \BFX_k\|$, we conclude that $M_k^v(\BFx)$ is a poor estimate of variance at $\BFx$.
That is, if $M_k^v(\BFx) \ge \sigmahat^2(\BFX_k,N_{k-1}) + c_v\Delta_k$, for some constant $c_v$ that depends on the unknown variance Lipschitz constant (set to an arbitrary value in implementation), then we will resort to $N_{k,1}(\BFx)=\lambda_k$. To see why, note that we can write 
\begin{align*}
    |M_k^v(\BFx)-\sigmahat^2(\BFX_k, N_{k-1})|&\leq|M_k^v(\BFx)-\sigma^2(\BFx)|+|\sigma^2(\BFx)-\sigma^2(\BFX_k)|+|\sigma^2(\BFX_k)-\sigmahat^2(\BFX_k, N_{k-1})|.
\end{align*} All the three terms here can be bounded by a factor of the trust-region size $\Delta_k$, the first term by a Taylor expansion, the second term due to the Lipschitz continuity assumption, and the third term due to the satisfied adaptive sample size at $\BFX_k$. 
This hybrid two-stage estimation algorithm is outlined in \Cref{alg:two-hybrid-new}.  



\begin{algorithm}[htp]  
\caption{$N_k(\BFx)$=\texttt{TwoStageEstimation($\Delta_k,\BFx,M_k^v,\kappa, L_v, \mcF_k$)}}
\begin{algorithmic}[1]
\label{alg:two-hybrid-new}
\REQUIRE trust-region radius $\Delta_k$, design point $\BFx$, variance model $M_k^v$, constant $c_v>0$, and history $\mcF_k$.
\IF{$M_k^v(\BFx) \ge \sigmahat(\BFX_k,N_{k-1}) + L_v\Delta_k $}
\STATE Call \Cref{alg:two-lambda}.
\ELSE
\STATE Call \Cref{alg:two-variance-mod}.
\ENDIF

\end{algorithmic}
\end{algorithm}

It is important to note that Algorithms \ref{alg:two-lambda}-\ref{alg:two-hybrid-new} are exclusively employed for new design points $\BFx$ that have not been previously evaluated. 
When $\BFx$ is a design point that has been previously evaluated, a variance estimate and initial sample sizes at $\BFx$ are already known.
Whenever we reevaluate points, we obtain the second-stage sample size directly as 
\begin{equation} \label{eq:ss-old}
    N_{k}(\BFx) =  \max\left\{N_{k-1}(\BFx),\lambda_k, \lambda_k\frac{\sigmahat^2(\BFx, N_{k-1}(\BFx))}{\kappa \Delta_k^4}\right\},
\end{equation}
where $N_{k-1}(\BFx)$ denotes the total number of simulation oracle calls performed up until the $(k-1)$th iteration. 

We remark that all our proposed two-stage sample sizes are subject to a probability of not actually satisfying the criterion in \eqref{eq:as2}. 
Therefore, they are theoretically suboptimal compared to employing \eqref{eq:sample_condition} to find the right sample size. 
Despite sacrificing the theoretical guarantees that come with determining an optimal amount of simulation effort, these two-stage sample sizes come with a pratical benefit of substantially reducing communications between the classical and quantum computers.

\section{Variance Model Informed-2STRO-DF (VMI-2STRO-DF)}
In this section, we will present the algorithm VMI-2STRO-DF, which incorporates the novel algorithmic components discussed 
in the two previous sections, namely
the variance model-informed point selection strategy (\Cref{sec:point-selection}) and the two-stage estimation processes (\Cref{sec:how-ss}). 
 Psuedocode for VMI-2STRO-DF is presented in \Cref{alg:VMI-2STRO-DF}. 
 Note that, for notational brevity, the index $m$ distinguishing VMI-2STRO-DF-$m$ correspond to which of the three distinct two-stage estimation algorithms is employed in \Cref{alg:VMI-2STRO-DF}.  

\begin{algorithm}[htp]  
\caption{\texttt{VMI-2STRO-DF-$m$}}
\label{alg:VMI-2STRO-DF}
\begin{algorithmic}[1]
\REQUIRE Initial incumbent $\BFx_{0}\in\mbR^d$, initial and maximum trust-region radius $\Delta_{0},\Delta_{\max}>0$, model fitness thresholds $0<\eta_1<\eta_2<1$ and certification threshold $\mu>0$, sufficient reduction constant $\theta>0$, expansion and shrinkage constants $\gamma_1>1$ and $\gamma_2\in(0,1)$, sample size lower bound sequence $\{\lambda_k\}$, adaptive sampling constant $\kappa>0$, and a Lipschitz constant estimate for the variance function $L_v$.
\FOR{$k=0,1,2,\hdots$}

\STATE \label{ASTRO:set-selection} \textit{Design Set Selection:} Select $\mcX_{k}=\{ \BFX_{k}^{i}\}_{i=0}^{2d}\subset\mcB(\BFX_{k};\Delta_{k})$ by calling \Cref{alg:var-model-const}.


\STATE \label{ASTRO:TRsubprob} \textit{Objective Model Construction:} For all $i=0,1,\dots, 2d$, estimate $\Fbar(\BFX_{k}^{i},N(\BFX_{k}^{i}))$ using the sample size determined by Algorithm $m$ and construct the model $M_{k}(\BFX)$ via \eqref{eq:syslineq}. 

\STATE \label{ASTRO:solveTRsubprob} \textit{Subproblem:} Approximately compute the trust-region model minimizer $\BFXtilde_{k+1}=\argmin_{\left\| \BFX-\BFX_k \right\| \leq \Delta_{k}}M_{k}(\BFX)$. 

\STATE \label{ASTRO:evaluate}\textit{Candidate Evaluation:} Define $\Ntilde_{k+1} = N(\BFXtilde_{k+1})$ and estimate $\Fbar(\BFXtilde_{k+1},\Ntilde_{k+1})$ using the sample size determined by Algorithm $m$. Define the best design point $\BFXhat_{k+1}=\argmin_{\BFx\in\mathcal{X}_{k}\backslash\BFX_k}\Fbar(\BFx, N_k(\BFx))$, its sample size $\hat{N}_{k+1} = N(\BFXhat_{k+1})$, incumbent's sample size $\Nhat_k = N\left(\BFX_k\right)$, direct search reduction $\Rhat_k=\Fbar(\BFX_{k},\Nhat_{k})-\Fbar (\BFXhat_{k+1},\Nhat_{k+1})$, subproblem reduction   $\Rtilde_k=\Fbar(\BFX_{k},\Nhat_{k})-\Fbar (\BFXtilde_{k+1},\Ntilde_{k+1})$, and model reduction $R_k=M_{k}(\BFX_{k})-M_{k}(\BFXtilde_{k+1})$. 
\STATE \label{ASTRO:ratio} \textit{Update:} Set
$(\BFX_{k+1},N_{k+1},\Delta_{k+1})=$ \[
\begin{cases}
    (\BFXhat_{k+1},\Nhat_{k+1},\gamma_1\Delta_{k}\land \Delta_{\max})& \text{if } \Rhat_{k}>\max\{\Rtilde_k,\theta\Delta_{k}^2\}, \\
    (\BFXtilde_{k+1},\Ntilde_{k+1},\gamma_1\Delta_{k}\land \Delta_{\max})              & \text{if }\Rtilde_{k}\geq\eta_2 R_{k} \text{ and }\mu\|\TD M_{k}(\BFX_{k})\|\ge\Delta_{k},\\
    (\BFXtilde_{k+1},\Ntilde_{k+1},\Delta_k)              & \text{if }\Rtilde_{k}\geq\eta_1 R_{k} \text{ and }\mu\|\TD M_{k}(\BFX_{k})\|\ge\Delta_{k},\\
    (\BFX_{k},\Nhat_{k},\gamma_2\Delta_{k})              & \text{otherwise},
\end{cases}
\] and $k=k+1$.
\ENDFOR
\end{algorithmic}
\end{algorithm}

\subsection{Construction of the local model $M_k^v$}
In each iteration of VMI-2STRO-DF, the construction of $M_k^v$ occurs prior to the selection of $\mcX$, implying that $M_k^v$ is always constructed using previously evaluated points, which is involved in the filtration $\mcF_{k-1}$. 
Hence, the deterministic coordinate basis geometry employed for constructing $\mcX_k$ cannot be used for $\mcX_k^v$. Moreover, especially at the beginning of a run of VMI-2STRO-DF, there may be an insufficient number of previously evaluated points to construct a interpolation or regression model $M_k^v$ that provides predictive accuracy on a current trust region. To address this gap, we propose a consistent approach for constructing the local model $M_k^v(\BFs)$ in \Cref{alg:var-model-const}. 
\Cref{alg:var-model-const} expands the trust region until the number of previously evaluated design points within the trust region exceeds $2d$. 
This condition is satisfiable in every iteration except the $0$th iteration.
We construct the variance model $M_k^v$ using all the design points within the expanded trust region $\mcX_k^v$, 
provided $\sfM(\Phi, \mcX_k^v)$ in \Cref{defn:polyintermd}, with the monomial basis $\Phi$ employed in \Cref{defn:diagHess} has a defined pseudoinverse. 
If $|\mcX_k^v|=2d+1$, we build a stochastic quadratic interpolation model with a diagonal Hessian. Otherwise, we construct a regression model with the same monomial basis employed in \Cref{defn:diagHess}. 
In instances where it is not possible to construct the variance model $M_k^v$, that is, in the first iteration or when $\sfM(\Phi, \mcX_k^v)$ does not admit a pseudoinverse, we elect not to construct $M_k^v$ at all. 
In such cases, we select the design set for $M_k$ using \Cref{alg:design-set}, the same approach used in history-informed ASTRO-DF for the purpose of model construction.
It is worth noting that the case where $\sfM(\Phi, \mcX_k)$ was not psuedoinvertible very rarely occurred after the $0$th iteration in our numerical experiments.


\begin{algorithm}[htp]  
\caption{[$\mcX_k$] = \texttt{VMI-ChooseDesignSet}($\BFX_k,\Delta_k, \mcF_k, w, k$)}
\label{alg:var-model-const}
\begin{algorithmic}[1]
\REQUIRE current iterate $\BFX_k$, trust-region radius $\Delta_k$, history $\mcF_k$, and some constant $w > 1$.
\STATE Select $\mcX_{k}=\{ \BFX_{k}^{i}\}_{i=0}^{2d}\subset\mcB(\BFX_{k};\Delta_{k})$ by calling \Cref{alg:design-set}.
\IF{$k > 0$}
\STATE Initialize $j=0$
\REPEAT  
\STATE Find the design set $\mcX^v_k$, which include all design points within the region $\mcB(\BFX_k,\widetilde{\Delta}_k^{(j)})$, where $\widetilde{\Delta}_k^{(j)} = \Delta_k w^j$.
\STATE Set $j = j+1.$
\UNTIL{$|\mcX^v_k| \ge 2d+1$.} 

\IF{$\sfM(\Phi, \mcX_k^v)$ is pseudoinvertible}
\STATE \textit{Variance Model Construction and Subproblem Solution:} Construct the stochastic quadratic regression/interpolation model with diagonal Hessian $M^v_{k}\left(\BFx\right)$ and compute an approximate minimizer $\BFXtilde_k^v := \argmin_{\left\| \BFx - \BFX_k\right\| \leq \Delta_{k}}M_k^v(\BFx)$.

\STATE \textit{Two-Stage Estimation:} Define $\widetilde{N}_k^v = N(\BFXtilde_k^v)$ and estimate $\Fbar(\BFXtilde_k^v, \widetilde{N}_k^v)$ by calling Algorithm $m$.

\STATE \textit{Design Set Update:} Find the closest design point to $\BFXtilde_k^v$ in the design set $\mcX_k$, i.e., $\widetilde{i} = \argmin_{i\in\{0,\dots,2d\}}\|\BFX_k^{i}-\BFXtilde_k^v\|$ and set $(\BFX_k^{\widetilde{i}},\Fbar(\BFX_k^{\widetilde{i}},N(\BFX_k^{\widetilde{i}}))=$ \[
\begin{cases}
    (\BFX_k^v,\Fbar(\BFXtilde_k^v, \widetilde{N}_k^v))& \text{if } \widetilde{i}\in\{1,2,3,\dots,2d\} and \mcX_k = \mcX_k^{cb}, \\
    (\BFX_k^v,\Fbar(\BFXtilde_k^v, \widetilde{N}_k^v))              & \text{else if }\widetilde{i}\in\{2,3,\dots,2d\}.\\
\end{cases}
\] 
\label{step:design-set}

\ENDIF

\ENDIF
\RETURN $\mcX_k$
\end{algorithmic}
\end{algorithm}

\subsection{Selection of the design set $\mcX_k$}
Similar to history-informed ASTRO-DF, VMI-2STRO-DF constructs a stochastic quadratic interpolation model with a diagonal Hessian, as in \Cref{defn:diagHess}.
This construction requires $2d+1$ design points. 
\Cref{alg:var-model-const} outlines the procedure for computing $\mcX_k$
Motivated by a desire to reuse some previously evaluated design points, but still maintain a well-poised geometry, 
\Cref{alg:var-model-const} begins with a scheme employed in history-informed ASTRO-DF, summarized in \Cref{alg:design-set}.
\Cref{alg:design-set} idenitifes all previously evaluated design points within the current trust region distinct from $\BFX_k$, and chooses one such design point, $\BFX_k^{1}$ maximally far from $\BFX_k$, breaking ties arbitrarily. 
\Cref{alg:design-set} then generates a basis of orthonormal vectors for $\mbR^d$ and completes the design set $\mcX_k$ by moving $\Delta_k$ in the positive and negative direction away from $\BFX_k$ in each orthonormal direction. 

Upon completion of \Cref{alg:design-set} in \Cref{alg:var-model-const}, $\mcX_k$ will generally contain $2d-1$ points not previously evaluated. 
Further motivated by keeping the total number of function evaluations low, 
after \Cref{alg:var-model-const} has obtained (and evaluated) a variance model minimizer $\BFXtilde_k^v$, we would like to include $\BFXtilde_k^v$ in $\mcX_k$. 
We replace the point in $\mcX_k$ closest to $\BFXtilde_k^v$ with $\BFXtilde_k^v$, provided  the selected point in $\mcX_k$ is not $\BFX_k^{1}$, which has already been evaluated (See \Cref{alg:design-set} in the Appendix). 
Therefore, at the termination of \Cref{alg:var-model-const}, the set $\mcX_k$ only contains $2d-2$ previously unevaluated points in the best case. 

\section{Numerical Results}

We will now assess and compare VMI-2STRO-DF with SO solvers. In this section, we will refer to history-informed ASTRO-DF as ASTRO-DF. The evaluation of SO solvers will be conducted using SimOpt \citep{eckman2022simopt}.
We implement various QAOA circuits using Qiskit \citep{Qiskit}. 

It is worth noting that in order to facilitate an effective comparison between the solvers, we have employed common random numbers (CRN) through SimOpt alongside the quantum simulator available in Qiskit. CRN is a method for variance reduction by querying the simulation-based oracle with the same random number stream. SimOpt offers versatile capabilities for applying CRN in diverse ways, enabling us to carry out sharper (with less variance) comparisons among the solvers. While CRN cannot be applied on a real quantum computer, our numerical experiments are conducted using the Qiskit quantum simulator, which does allow us to fix the random number seed.
We reiterate that nothing in VMI-2STRO-DF assumes access to CRN streams, and this choice in experimentation was made only to reduce sources of variance in comparing across solvers. 


\subsection{Comparison of the two-stage estimation algorithms}
Before comparing VMI-2STRO-DF with other SO solvers such as SPSA and Nelder-Mead, we perform a comparison of the three distrinct variants of two-stage estimation strategies Algorithms \ref{alg:two-lambda}-\ref{alg:two-hybrid-new}), as indexed by VMI-2STRO-DF-$m$ in \Cref{alg:VMI-2STRO-DF}. 

\begin{figure} [htp]
\centering
\subfloat[$c_n = 0$ and $c_s = 1$]{%
\resizebox*{8cm}{!}{\includegraphics{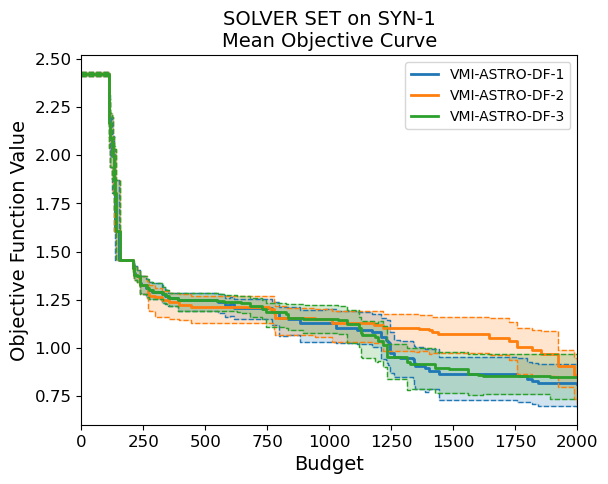}}\label{fig:himme-0}}
\subfloat[$c_n = 1000$ and $c_s = 1$]{%
\resizebox*{8cm}{!}{\includegraphics{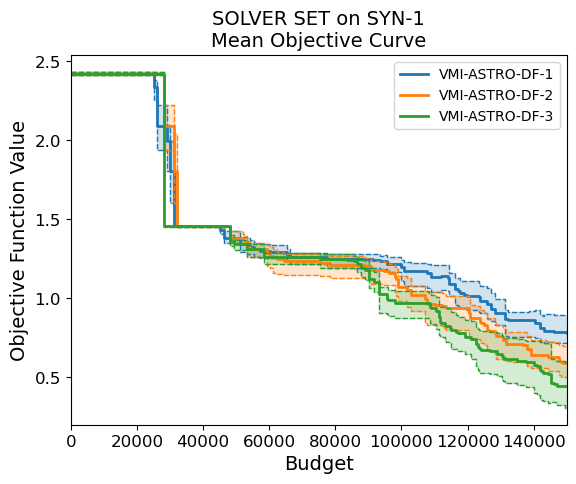}}\label{fig:himme-1000}} 
\caption{Performance of VMI-2STRO-DF-$m$ for $m\in\{1,2,3\}$ on the stochastic Himmelblau function.
Solid lines denote mean performance over 20 macro-replications, bands denote 95\% confidence intervals.
Initial incumbent is $\BFX_0=(-5,-5)$. 
The $x$-axis shows the computational burden as measured in \eqref{eq:budget}
, and the $y$-axis shows expectation function value on a log scale.}
\label{fig:himme-three-two-stage}
\end{figure}

As a first test problem, we again consider the stochastic Himmelblau function described in \eqref{eq:himme-problem}, but with a higher variance of the stochastic noise $\xi~\mcN(0,10|(x_1-3)(x_2-2)|)$. 
To guarantee both robustness and accuracy, we ran each algorithm 20 times. 
Based on our preliminary experiment in \Cref{sec:point-selection}, we expect that VMI-2STRO-DFs will eventually converge to the global optimum, provided the two-stage estimation algorithm employed provides a reasonable estimate of variance. In the numerical experiments, the metric for the computational budget takes into account both the costs ($c_s$) related to acquiring a single sample and the communication costs ($c_n$). To calculate the total computational budget expended during the experiments, we use the  formula
\begin{equation} \label{eq:budget}
    c_n Q_n + c_s W_s,
\end{equation}
where $W_s$ quantifies the total count of oracle calls (shots) requested by a SO method, where in our context $W_s=\sum_{k=1}^{T}\sum_{i=0}^p N_k^p$ with $T$ being the number of iterations, and $Q_n$ quantifies the total count of communications made with the quantum computer.
See \Cref{fig:himme-three-two-stage}.
When the communication cost $c_n$ is assumed low relative to $c_s$, the choice of two-stage sampling strategy appears to make little difference. 
However, as communication costs increase,  VMI-2STRO-DF-3, which employs the hybrid two-stage sampling strategy, demonstrates superior performance as depicted in \Cref{fig:himme-1000}. 

\begin{figure} [htp]
\centering
\includegraphics[width=0.4\columnwidth]{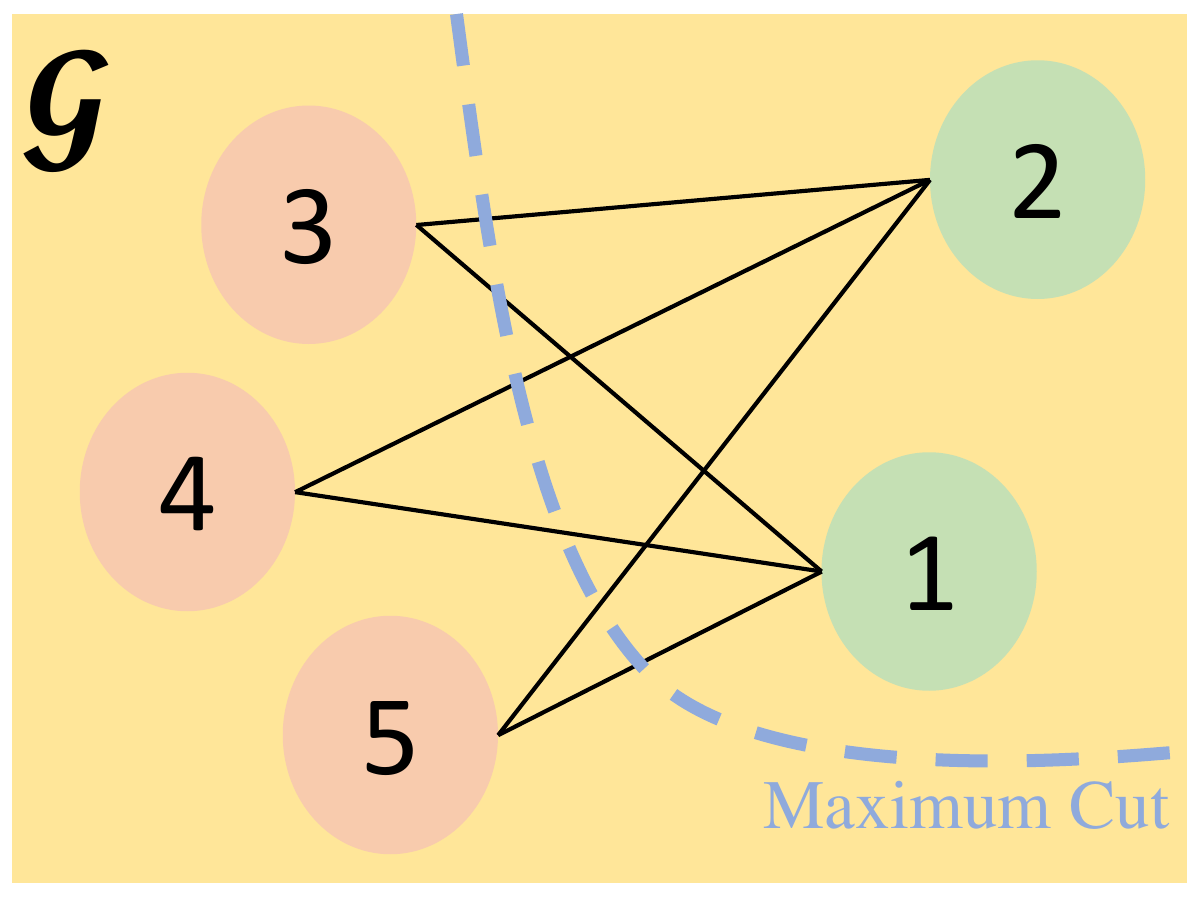}
\caption{An illustration of the graph of the Max-Cut problem. We aim to discover a way to partition the vertices of the graph into two complementary sets in such a manner that maximizes the number of edges between these two sets.} 
\label{fig:maxcut-problem}
\end{figure}

Our second test is on a QAOA circuit for solving a Maxcut problem \citep{farhi2014qaoa}. Specifically, we have implemented the standard QAOA with a depth of 5. 
We now provide a brief description of the Maxcut problem and implementation. Let us consider a graph $\mcG$ represented as $\mcG = [\mcV,\mcD]$, where $\mcV$ denotes the set of vertices and $\mcD$ represents the set of edges. 
A \textit{cut} of the graph $\mcG$ is defined as a partition $(\mcV_1, \mcV_2)$ of the graph vertices such a way that there are no edges between $\mcV_1$ and $\mcV_2$.
A \textit{maximum cut} of a graph $\mcG$ is a cut such that, among all possible cuts, the number of edges between $\mcV_1$ and $\mcV_2$ is maximized. 
This problem can be formulated as a quadratic unconstrained binary optimization (QUBO) problem, 
\begin{equation}
   \max \frac{1}{2}\sum_{i,j \in \mcV}(1-x_ix_j),
\end{equation} \label{eq:maxcut-origin}
where $\mcV = \{1,2,\dots,5\}$ and $x_i \in \{-1,+1\}$ for $i \in \mcV$.
\Cref{fig:maxcut-problem} provides a visual representation of a maximum cut on a toy graph. 

\begin{figure} [htp]
\centering
\subfloat[$c_n = 0$ and $c_s = 1$]{%
\resizebox*{8cm}{!}{\includegraphics{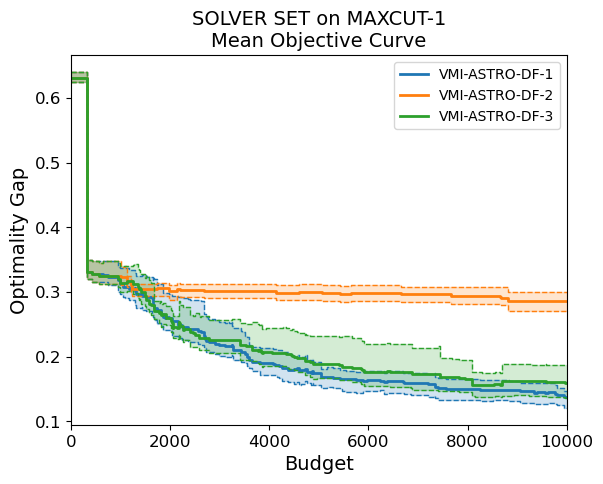}}\label{fig:maxcut-0}}
\subfloat[$c_n = 10$ and $c_s = 1$]{%
\resizebox*{8cm}{!}{\includegraphics{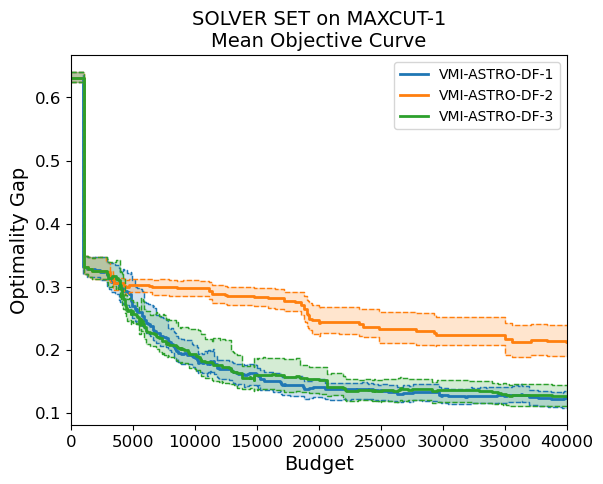}}\label{fig:maxcut-10}}

\subfloat[$c_n = 100$ and $c_s = 1$]{%
\resizebox*{8cm}{!}{\includegraphics{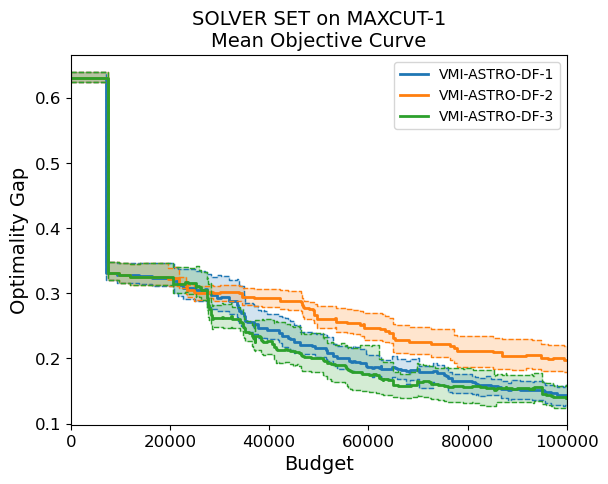}}\label{fig:maxcut-100}} 
\subfloat[$c_n = 100$ and $c_s = 1$]{%
\resizebox*{8cm}{!}{\includegraphics{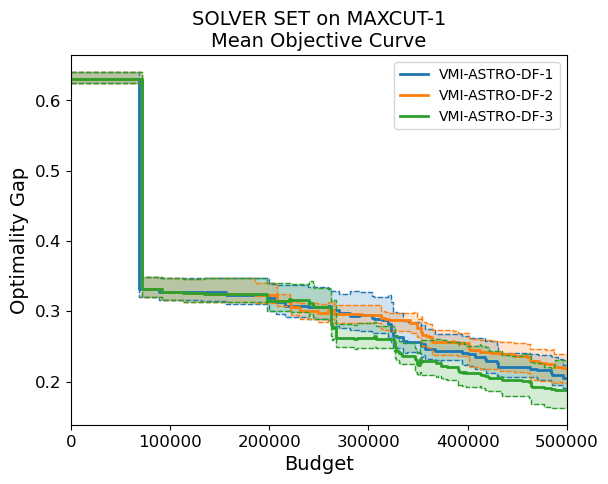}}\label{fig:maxcut-100}} 
\caption{Finite time performance on Max-Cut problem with 95\% confidence interval. The $x$-axis shows the computational burden as measured in \eqref{eq:budget}. Y-axis shows the optimality gap on the log scale.}
\label{fig:maxcut-three-two-stage}
\end{figure}

We conducted tests on three VMI-2STRO-DF variants for the Max-Cut problem, utilizing 20 macro-replications as illustrated in \Cref{fig:maxcut-three-two-stage}. Similar to the outcome observed with the stochastic Himmelblau function, VMI-2STRO-DF-3 exhibits superior performance as communication costs increase. Furthermore, as the communication costs increase, the performance gap between VMI-2STRO-DF-1 and VMI-2STRO-DF-2 narrows. 
It is worth noting that when communication costs are zero, VMI-2STRO-DF-1 can demonstrate strong performance with a careful selection of $\lambda_k$, see \Cref{fig:maxcut-0}. 
However, the performance of VMI-2STRO-DF-1 appears to be highly sensitive to various factors, including the initial incumbent $\BFX_0$ and choices of $\{\lambda_k\}$.

\subsection{Comparison with stochastic optimization solvers}
\begin{figure} [htp]
\centering
\subfloat[$c_n = 0$ and $c_s = 1$]{%
\resizebox*{8cm}{!}{\includegraphics{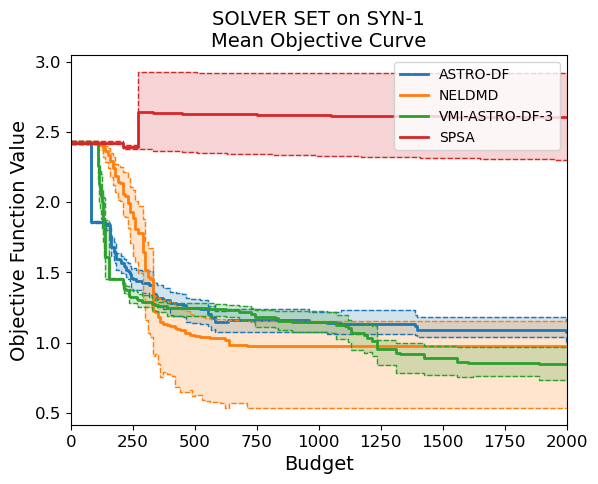}}\label{fig:himme-0-entire}}
\subfloat[$c_n = 1000$ and $c_s = 1$]{%
\resizebox*{8cm}{!}{\includegraphics{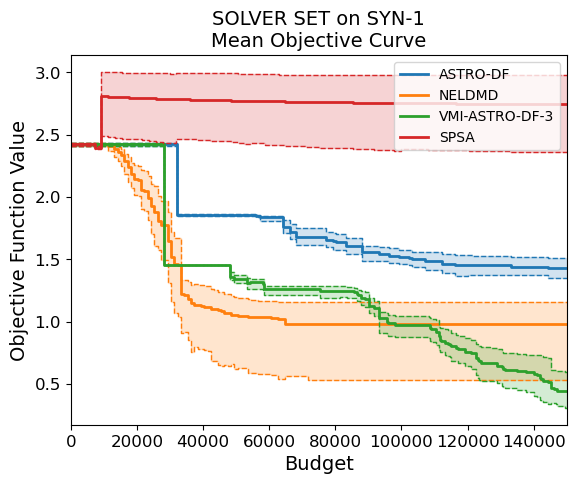}}\label{fig:himme-1000-entire}} 

\caption{Finite time performance on the stochastic Himmelblau function with 95\% confidence interval with initial design point = (-5,-5). X-axis shows the computational burden with $c_n$ and $c_s$, where $c_n$ is the communication costs and $c_s$ is the costs associates with obtaining single sample. Y-axis shows the objective function value on the log scale.}
\label{fig:himme-entire}
\end{figure}

\begin{figure} [htp]
\centering
\subfloat[$c_n = 0$ and $c_s = 1$]{%
\resizebox*{8cm}{!}{\includegraphics{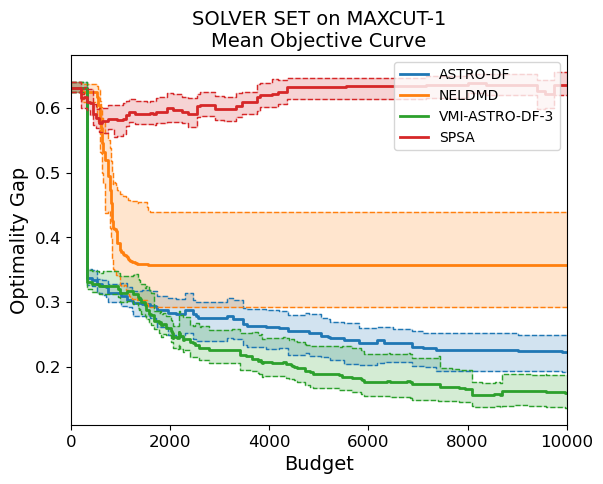}}\label{fig:maxcut-0-entire}}
\subfloat[$c_n = 10$ and $c_s = 1$]{%
\resizebox*{8cm}{!}{\includegraphics{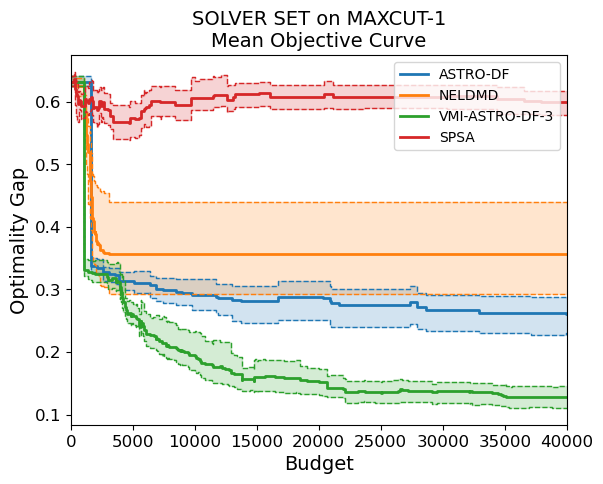}}\label{fig:maxcut-10-entire}}\\

\subfloat[$c_n = 100$ and $c_s = 1$]{%
\resizebox*{8cm}{!}{\includegraphics{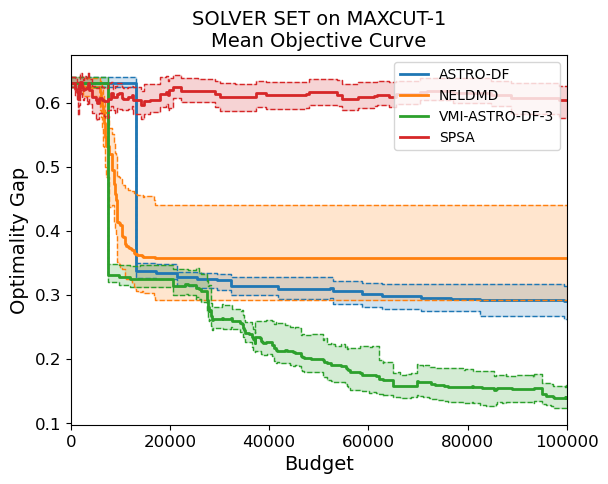}}\label{fig:maxcut-100-entire}}
\subfloat[$c_n = 1000$ and $c_s = 1$]{%
\resizebox*{8cm}{!}{\includegraphics{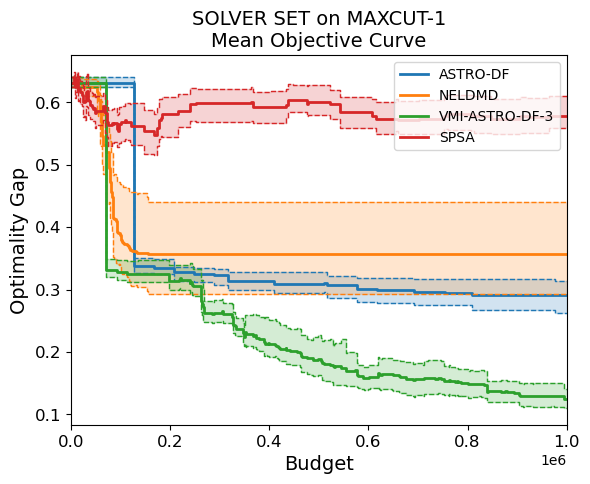}}\label{fig:maxcut-1000-entire}}\\
\caption{Finite time performance on Max-Cut problem with 95\% confidence interval. X-axis shows the computational burden with $c_n$ and $c_s$, where $c_n$ is the communication costs and $c_s$ is the costs associates with obtaining single sample. Y-axis shows the optimality gap on the log scale.}
\label{fig:maxcut-entire}
\end{figure}

In this section, we will discuss the comparison between VMI-2STRO-DF-3 with ASTRO-DF, SPSA, and Nelder-Mead on two previous problems. On both problems, VMI-2STRO-DF-3 consistently outperforms other solvers in terms of finding the optimal solution, regardless of communication costs. This superiority is evident in \Cref{fig:himme-entire} and \Cref{fig:maxcut-entire}. Even when communication costs are zero, VMI-2STRO-DF-3 converges to a superior solution more quickly than its competitors. This remarkable performance can be attributed to the unique combination of a design set that incorporates the minimizer of the variance model and the direct search method, as elaborated in \Cref{sec:point-selection}. Furthermore, as communication costs escalate, the performance gap between VMI-2STRO-DF-3 and other solvers widens. More precisely, Nelder-Mead and SPSA tend to become stuck at solutions significantly distant from the optimal one, and ASTRO-DF is limited in the number of iterations it can perform due to the substantial communication costs.  

\section{Conclusion}
In the optimization problem that forms an essential part of VQAs, two phenomena come to the forefront.
The first pertains to latency, which leads to an increase in the time required to acquire a single sample of shots. 
Given its substantial computational overhead during the optimization process, it is imperative for the optimizer to be thoughtfully designed with a focus on minimizing the quantum computer access frequency. 
The second phenomenon is the gradual reduction of true variance that accompanies the decaying optimality gap. 
In this paper, we introduced a novel stochastic trust region method to tackle VQA optimization problems, leveraging the distinctive characteristics of diminishing variance and communication costs.
We named this method VMI-2STRO-DF. 
To leverage the characteristic of diminishing variance, VMI-2STRO-DF constructs two local models using function estimates and variance estimates. The minimizer of the variance model is included in the design set, helping us to find better solutions at each iteration and the global optimum. 
Moreover, to reduce communication costs, VMI-2STRO-DF restricted the number of accesses to the quantum computer to two accesses per design point at each iteration. 
Our numerical results showcase the effectiveness of VMI-2STRO-DF on two problems from QAOA. Even in scenarios without communication costs, our approach outpaces other solvers in terms of finding quality solutions more efficiently.
Moreover, as communication costs increase, the performance gap between VMI-2STRO-DF and other solvers widens, further highlighting its advantage in latency-constrained settings.

\section*{Acknowledgments}
This material is based upon work supported by the U.S. Department of Energy, Office of Science, National Quantum Information Science Research Centers and 
the Office of Advanced Scientific Computing Research, Accelerated Research for Quantum Computing program under contract number DE-AC02-06CH11357. Sara Shashaani and Yunsoo Ha also gratefully acknowledge the U.S. National Science Foundation for support provided by grant number CMMI-2226347.  

\bibliography{main-paper} 

\newpage
\section*{Appdendix. Algorithm for selecting $\mcX_k$ in history-informed ASTRO-DF}
\label{appen:picking-algorithm}

\begin{figure} [htp]
\centering
\subfloat[iteration $k-1$]{%
\resizebox*{7cm}{!}{\includegraphics{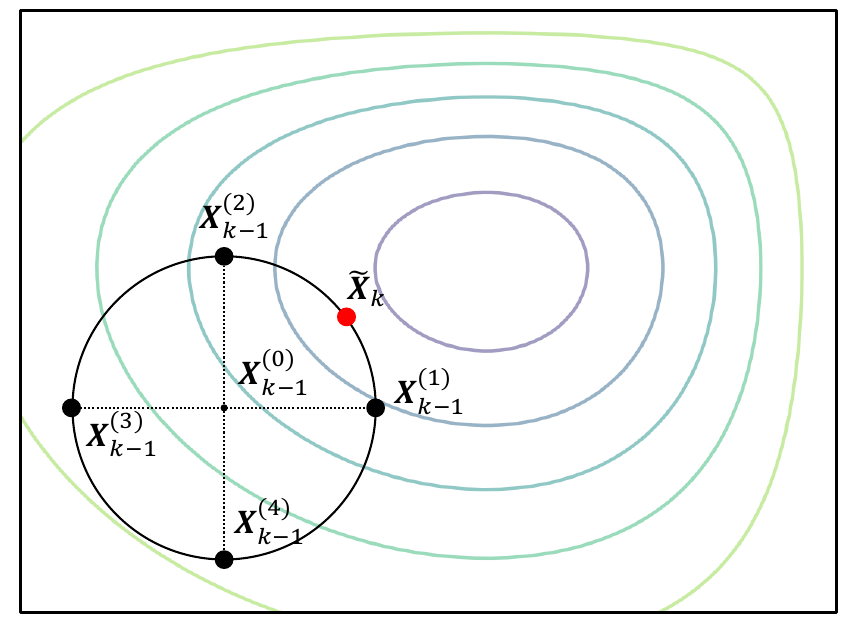}}\label{fig:2_A}}
\subfloat[iteration $k$]{%
\resizebox*{7cm}{!}{\includegraphics{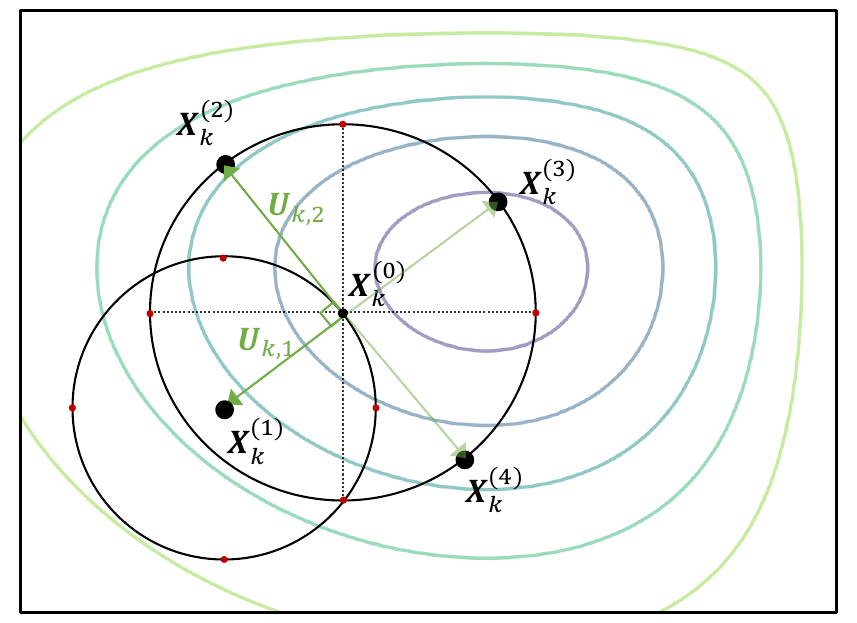}}\label{fig:2_B}}
\caption{An example demonstrating the application of the rotated coordinated basis \citep{ha2023wsc}. \Cref{fig:2_A} illustrates the coordinate basis using a coordinate system defined by elementary basis vectors, which is history-informed ASTRO-DF's default coordinate basis in the absence of reusable design points within the trust region. 
\Cref{fig:2_B} illustrates a rotated coordinate basis. 
In this case, in the $k$th iteration of history-informed ASTRO-DF, the design point $\BFX_{k-1}$ is the farthest from $\BFX_k$ among all previously evaluated design points, and so we choose $\BFX_k^{1} = \BFX_{k-1}$. Orthogonalizing against $\BFU_{k,1}:= \BFX_k - \BFX_k^{1}$ deterministically defines the rotated coordinate basis in the $k$th iteration.}
\label{fig:rcb} 
\end{figure}

\Cref{alg:design-set} shows the psuedocode for selecting the design set at iteration $k$ in history-informed ASTRO-DF \citep{ha2023wsc}. The example illustrating the application of \Cref{alg:design-set} is detailed in \Cref{fig:rcb}.

\begin{algorithm}[htp]  
\caption{$\mcX_k$=\texttt{ChooseDesignSet($\Delta_k,\BFX_k,\mcF_k$)}}
\label{alg:design-set}
\begin{algorithmic}[1]
\REQUIRE trust-region radius $\Delta_k$, iterate $\BFX_k$.
\STATE Find all previously evaluated design points within the trust region $\mcB(\BFX_k,\Delta_k)$, denote it $\mcR_k$.
\IF{$\mcR_k=\{\BFX_k\}$}
\STATE  Select the design set $\mcX_k=\mcX_k^{cb}$ following \Cref{defn:diagHess}.
\ELSE
\STATE Select $\BFX_k^{1}$ via
\begin{equation*}
    \BFX_k^{1}=\argmax_{\BFx \in \mcR_k} \|\BFX_k - \BFx\|_2 = \BFX_k + P_k \BFU_{k,1},
\end{equation*}
where $\|\BFU_{k,1}\| = 1$ and $P_k=\|\BFX_k^{1} - \BFX_k^{0}\|$.
\STATE Compute a set of $d-1$ vectors $\{\BFU_{k,2}, \dots, \BFU_{k,d}\}$ mutually orthonormal to $\BFU_{k,1}$ 
\STATE $\mcX_k \gets \{\BFX_k, \BFX_k + P_k\BFU_{k,1}, \BFX_k + \Delta_k\BFU_{k,2}, \dots, \BFX_k + \Delta_k\BFU_{k,d}, \BFX_k - \Delta_k\BFU_{k,1}, \dots, \BFX_k - \Delta_k\BFU_{k,d} \}$.
\ENDIF
\STATE \textbf{Return} $\mcX_k$.
\end{algorithmic}
\end{algorithm}

\end{document}